\documentclass[12pt,leqno,fleqn,epsfig]{article}
\usepackage{amssymb, epsfig, amsmath, amsthm}
\usepackage{mathrsfs}
\usepackage{color}       

\newcommand{\R}{\mathbb{R}}
\newcommand{\PP}{\mathbb{P}}

\newcommand{\Z}{\mathbb{Z}}

\newcommand{\N}{\mathbb{N}}

\renewcommand{\a}{\mu}

\newcommand{\mes}{\operatorname{\rm meas}}    
\newcommand{\esssup}{\operatorname*{ess\,sup}}

\newcommand{\supp}{\operatorname*{supp}}

\newcommand{\bb}{\begin{equation}}
\newcommand{\ee}{\end{equation}}
\newcommand{\bq}{\begin{eqnarray}}
\newcommand{\eq}{\end{eqnarray}}
\newcommand{\bqn}{\begin{eqnarray*}}
\newcommand{\eqn}{\end{eqnarray*}}

\newcommand{\var}{\varepsilon}

\newcommand{\intl}{\int\limits}

\newcommand{\Beweisende}{\rule{0.2cm}{0.2cm}}

\newcommand{\intmw}{{\int\hspace{-830000sp}-\!\!}}

\newcounter{secnum}

\newtheorem{thm}{Theorem}[section]
\newtheorem{cor}[thm]{Corollary}
\newtheorem{lem}[thm]{Lemma}

\theoremstyle{definition}

\newtheorem{defin}[thm]{Definition}
\newtheorem{rem}[thm]{Remark}

\title{Energy concentrations  and  Type I blow-up \\for  the 3D Euler equations} 
 
\author{Dongho Chae$^*$  and J\"{o}rg Wolf $^\dagger$\\
\ \\
Department of Mathematics\\
Chung-Ang University\\
 Seoul 156-756, Republic of Korea\\
 ($*$)e-mail: dchae@cau.ac.kr\\
($\dagger$)e-mail: jwolf2603@cau.ac.kr}
\date{}
\begin{document}
\maketitle
\begin{abstract}
We exclude Type I blow-up, which occurs in the form of  atomic concentrations of the $L^2$ norm for the solution of  the 3D incompressible  Euler equations.  As a corollary  we prove nonexistence of discretely self-similar blow-up  in the energy conserving scale. 
\\
\
  \noindent{\bf
keywords:} incompressible Euler equations, finite time blow-up, energy concentration, discretely self-similar solution

\end{abstract}

\section{Introduction}
\label{sec:-1}
\setcounter{secnum}{\value{section} \setcounter{equation}{0}
\renewcommand{\theequation}{\mbox{\arabic{secnum}.\arabic{equation}}}}

We consider the $ n$-dimensional  Euler equations   in  $ \R^{n}\times (0, +\infty)$
\begin{equation}
\left\{\aligned &\partial _t v  + (v \cdot \nabla ) v  = - \nabla p,\quad  \nabla \cdot  v =0,\\
 &\quad v(x,0)=v_0(x),
 \endaligned 
 \right.
\label{1.1}
\end{equation}
where $v=(v_1 (x,t), \ldots, v_n (x,t))$,  $(x,t)\in  \R^n \times (0, +\infty)$.
For the Cauchy problem of the system \eqref{1.1} the local well-posedness 
 in the setting of standard Sobolev space $H^m (\Bbb \R^n )$, $ m> \frac{n}{2} +1$,  is proved by Kato in \cite{kat}.  The question of finite time blow-up of such local in time classical solution, however, is an outstanding open problem in the mathematical fluid mechanics(see e.g. \cite{maj, con1} for an introduction and surveys of partial results on the  problem, and \cite{gra, gre, ker, luo} for the related numerical works).  In this direction of study there are also well-known results on the blow-up criterion\cite{bea, con2, den, koz}, where the authors deduced various sufficient conditions for the blow-up. We also mention a recent result
 by Tao\cite{tao}, which shows  the blow-up for 
 a model equation having   similar conservation properties to the Euler system.   \\

The  aim of the present paper is to study the possibility of the finite time blow-up in terms of the energy concentrations in  the 3D Euler equations.  The  phenomena of $L^2$ norm concentration at the  blow-up time is well-known in the other nonlinear evolution equations. For example in the nonlinear Schr\"{o}dinger equation it is found that there exists  a solution  which shows that the mass($L^2$ norm of the solution) is concentrating in the form of finite sum of Dirac measures at the  blow-up time\cite{mer1, mer2}. Similarly, in the chemotaxis equation  the $L^1$ norm of solution is shown to be evolved into Dirac measures  in the finite time for a sufficiently large initial data \cite{her}.  We also find that there exists a study of the energy concentration for the Navier-Stokes equations, in the context different from ours  in \cite{arn}.\\

In our case of the 3D Euler system, under Type I condition for  the velocity gradient we are able to exclude the atomic concentrations of velocity $L^2$ norm  at the possible blow-up time. This means that there exists no concentration of the energy into  isolated points in $\Bbb R^n$ at the possible blow-up time if we assume
Type I condition for the blow-up rate. As an immediate corollary of this result we exclude the discretely self-similar(DSS) blow-up in the energy conserving scale. 

Let us denote by $L^2_\sigma (\Bbb \R^n )$ the closure of $\{ \varphi \in C_c ^\infty (\Bbb \R^n )\, | \, \nabla \cdot \varphi =0\}$ in $L^2 (\Bbb \R^n )$. 
Given a domain $\Omega \subset \Bbb R^n$,  we denote by $ \mathcal{M}(\Omega )$ the space of all bounded Radon measures $ \mu \in C^0 _c(\Omega )^{ \ast}$. 
The space $ \mathcal{M}(\Omega )$ will be equipped with the norm 
\[
\| \mu \|_{ \mathcal{M}} =  \sup_{ \substack{\varphi \in C^0_c(\Omega )\\ \max_{ \Omega } | \varphi | \le 1 }} 
\intl_{\Omega } \varphi d\mu.   
\]
In particular, by $ \mathcal{M}^+(\Omega )$ we denote the subspace of all nonnegative $ \mu \in \mathcal{M}(\Omega )$, i.e. 
\[
\intl_{\Omega } \phi d\mu \ge 0\quad  \forall\,\phi \in C^0_c(\Omega ), \quad  \text{ with}\quad  \phi \ge 0. 
\]
If $ f\in L^\infty(a,b; L^1(\Omega ))$, $-\infty < a<b < +\infty$,  by $ \mathcal{M}_f(b)$ we denote the set of all  
$ \sigma_0  \in \mathcal{M}(\Omega )$ 
such that  there exists a sequence $ \{ s_k\} $ in the set of Lebesgue  points of $ f(\cdot )$ such that  $ s_k \rightarrow b$ as 
$ k \rightarrow +\infty$ and 
\[
f(s_k)  d x \rightarrow 
\sigma _0 \quad \text{weakly-$  \ast$ in }\quad  \mathcal{M}(\Omega  )\quad   \text{ as}  \quad  k\rightarrow +\infty.
\]
Here $ t\in (a,b]$ is called  a Lebesgue point of $ f(\cdot )$ if 
\[
\frac{1}{h}  \intl_{t-h}^{t} f(s)  ds \rightarrow f(t)\quad  \text{ {\it in}}\quad L^1(\Omega )\quad  \text{ {\it as}}\quad 
h \rightarrow 0^+. 
\]
 Note that due to Lebesgue's differentiation theorem for the Bochner integrable functions (see  e.g. 
  in \cite[Theorem 2, pp.\,134]{boc}) almost every $ t\in (a,b]$ is a Lebesgue  point of $ f(\cdot )$. 
 
\vspace{0.2cm} 
For simplicity  in the discussion  below we consider our time interval $(-1,0)$, and  fix $t=0$ as the possible blow-up time.
Our main theorem of this paper is the following.

\begin{thm}
\label{thm1.4} 
Let $ \Omega \subset \R^n $ be a domain. 
Let $ v\in L^\infty(-1,0; L^2_{ \sigma }(\Omega )\cap L^\infty_{ loc}([-1,0), W^{1,\, \infty}(\Omega ))$ be a solution  to 
\eqref{1.1} in $ \Omega \times (-1,0)$ satisfying the following Type I blow-up condition at $ t=0$
\begin{equation} 
\label{rate}  \qquad \sup_{ t\in (-1, 0)} (-t) \| \nabla v(t)\|_{ L^\infty(\Omega )} < +\infty.
\end{equation}

Then every measure $ \sigma_0  \in \mathcal{M}_{ | v|^2}(0)$ 
 has  no atoms, i.e. 
\begin{equation}
\sigma_0  (\{ x\})=0\quad  \forall\,x \in \Omega.
\label{noAtom}
\end{equation}  
If in addition, $ v(t) \rightarrow v_0$ weakly in $ L^2(\Omega )$ as 
$ t \rightarrow 0^-$ for some $ v_0\in L^2(\Omega )$, then $\mathcal{M}_{ | v|^2}(0)= \{ \sigma _0 \}$,
\begin{equation}
| v(t)|^2 dx \rightarrow 
\sigma _0  \quad \text{weakly-$  \ast$ in }\quad  \mathcal{M}(\Omega  )\quad   \text{ as}  \quad  t \rightarrow  0, 
\label{convm}
\end{equation}
and $ \sigma _0 $ has no atoms. 

\end{thm}

 In the case $\Omega=\Bbb R^n$ in the above theorem  the fact $p\in L^{\frac32}_{\rm {loc}} (\Bbb R^n\times (-1, 0))$ follows from 
 by the Calder\'on-Zygmund inequality and the velocity-pressure relation, $\Delta p= -
 \sum_{i,k=1}^{n}\partial _j\partial _k (v_j v_k)$. Therefore, as an immediate consequence of  Lemma 2.3 (with $g=v, f=0$) below the set  $ \mathcal{M}_{ | v|^2}(0)$ contains only one element $ \sigma _0\in \mathcal{M}(\R^n )$, which gives the  following.
  
  \begin{cor}
\label{cor1.2}
Let $  v\in L^\infty(-1,0; L^2_{ \sigma }(\R^n ))\cap L^\infty_{ loc}([-1,0); W^{1,\, \infty}(\R^n ))$ be a solution of the Euler equations  \eqref{1.1}  satisfying  \eqref{rate}
with $\Omega = \Bbb R^n$.  Then, there exists $\sigma_0 \in \mathcal{M}(\R^n )$ such that 
\begin{equation}
| v(t)|^2 dx \rightarrow 
\sigma _0  \quad \text{weakly-$  \ast$ in }\quad  \mathcal{M}(\R^n )\quad   \text{ as}  \quad  t \rightarrow  0^-, 
\label{cor12}
\end{equation}
and 
$
\sigma_0  (\{ x\})=0
$ for all $x\in \R ^n$.
\end{cor}   
  
\begin{rem}
In particular, under Type I condition the limiting measure of the form $\sigma_0= \sum_{k=1} ^\infty c_k \delta_{x_k} + fdx$ 
with a sequence $\{c_k\}_{k=1}^\infty$   of nonnegative constants and $f\in L^1_{\rm{loc}} (\Bbb R^n)$,  is excluded  contrary to the case of the nonlinear Schr\"{o}dinger equation\cite{mer1, mer2} and the chemotaxis equation\cite{her}.  Currently, we are not able to exclude the possibility of energy concentration into a set of positive Hausdorff dimension under Type I condition, which would be an interesting subject for future study.

\end{rem}
\begin{rem}
\label{rem1.4}
We note that for $n=3$ we have
$$L^\infty (-1,0;  L^2 (\Bbb \R^3 ))\cap L^\infty_{\rm{loc}} ([-1,0); W^{1, \infty}(\Bbb \R^3 )) \subset L^3_{\rm{loc}}( [-1, 0); B^{\frac13 +s}_{3, \infty} (\Bbb R^3)), \quad \forall s \in  \Big(0, \frac23\Big]$$
which is the energy conserving class for the weak solutions to the Euler equations $v( \cdot, t)$ for $t\in [-1, 0)$ as  studied in \cite{con4}. As $t\to 0$, however, we cannot say anything about the energy conservation, and the existence of a definite particle trajectory map  up to $t=0$.  Therefore,  the energy concentration to a general measure zero set at the blow-up time cannot be excluded by  a naive application of the volume preserving property of the particle trajectory map.
\end{rem}
In order to discuss an implication of the above theorem on  the scenario of the discretely self-similar blow-up we first recall that  a solution $(v,p)$ of the Euler equations is self-similar if there exists $\alpha\not=-1$ such that 
\begin{equation} \label{ss}
v(x,t)=\lambda ^{\alpha} v(\lambda x, \lambda ^{\alpha+1} t), \quad  p(x,t)=\lambda ^{2\alpha}p(\lambda x, \lambda ^{\alpha+1} t)\quad \forall (x,t)\in \Bbb \R^n \times (0, +\infty)
\end{equation}
for  all $\lambda >1$.  The discrete self-similarity is a  more general concept; a solution $(v,p)$ of the Euler equations is called discretely self-similar(we say $\lambda$-DSS), if there exists $\alpha\not=-1$  and $\lambda >1$  such that \eqref{ss} holds. 
There have been previous studies on the exclusion of the
scenario of  self-similar blow-up \cite{cha1, cha2, cha3} in the Euler equations. Note that discretely self-similar solutions preserve the energy only if $\alpha=n/2$, which is called the energy conserving scale.
The  previous studies on the exclusion of discretely self-similar blow-up scenarios were mostly done in the other cases than the  energy conserving scale, for which the solution belongs to $L^q(\Bbb R^n)$, $q\not=2$, mainly due to the difficulties to prove Liouville type theorems for the corresponding profile equations.  The question of nonexistence of 
self-similar and/or discretely self-similar singularities in the 3D Euler equations has been open only in this case of energy conserving scale, while all the other cases are excluded under suitable decay conditions at infinity on the profiles\cite{cha1, cha2, cha3}.\\
 \\

As  proved in Section 5 below  the $\lambda$-DSS blow-up in the case $\alpha=n/2$ is a special case of the one point energy concentration at the time of blow-up.    
 Therefore as a consequence of Corollary\,\ref{cor1.2} we exclude the scenario of DSS blow-up in the energy conserving scale as follows.
 
\begin{cor}
\label{cor1.3}
Let $  v\in L^\infty(-1,0; L^2_{ \sigma }(\R^n ))\cap L^\infty_{ loc}([-1,0); W^{1,\, \infty}(\R^n ))$ be a solution of the Euler equations \eqref{1.1} satisfying  \eqref{rate}. 
If $ v$ is $ \lambda$-DSS solution with the energy conserving scale, i.e. if there exists $ \lambda >1$ such that  
\begin{equation} \label{rep}
v(x,t)=\lambda^{\frac{n}{2}} v(\lambda x, \lambda ^{ \frac{n+2}{2}} t) \quad  \forall\,(x,t)\in \R^n \times (-1,0).
\end{equation} 
Then $ v \equiv 0$.   

\end{cor} 

\hspace{0.5cm}
The paper is organized as follows. In Section\,2 we recall the notion of local pressure for bounded domains  and  exterior domains as well, which was previously introduced in \cite{wol} for the Navier-Stokes equations. 
Here the pressure gradient will be written as $ \nabla p = \partial _t \nabla p_h + \nabla p_0$ in the sense of distribution, where $ \nabla p_h$ stands for the harmonic part  associated to $ v$,  while $ p_0$ represents the  part associated to $ (v\cdot \nabla) v$. 
This eventually leads to a  local energy inequality in terms of the  new  localized energy  $  \intl_{\Omega } | \widetilde{v }(t) |^2 \phi dx$  with a cut-off function $\phi$,
where $ \widetilde{v} = v+ \nabla p_h$. A solution satisfying  this form of the local energy inequality  will be called {\em local suitable weak solution} as it has been introduced by 
Definition\,2.1.  As an important consequence   this  notion we show that for such solutions the energy $ | \widetilde{v} (t)|^2$ admits a unique measure valued     
trace,  which is in fact weak-$ \ast$  left continuous.   Furthermore,  we show that the question of concentration of $ | v(t)|^2$   at  the possible blow-up time can be reduced 
to that of  concentration of $ | \widetilde{v} (t) |^2$.  Section\,3 is devoted to  special case of  removing one point concentration of the energy 
for solution to the Euler equations in the whole space of $\Bbb R^n$. In particular, we are able to prove Theorem\,1.1 for this restricted situation,  which is stated in Thoerem\,\ref{thm3.1}. The proof  of Theorem 3.1 is based on several space-time  decay properties of the velocity field as $ t \rightarrow 0^-$. The proof of the decay estimates are  presented in Subsections 3.2 and 3.3. In particular, in Subsection\,3.3 we show that the energy $ | \widetilde{v}(t) |^2$  of any exterior 
subdomain excluding the concentration point converges to zero with   arbitrary  polynomial order.  In Subsection\,3.4  we complete the proof of Theorem\,3.1 based on a local estimate for the function $ w=v((-t)^{ \theta }x, t)$ for a suitable 
$ 0 < \theta < 1$,  which  by virtue of  Gronwall's lemma yields triviality of $ w$ in an exterior domain. Next, in Section\,4 
 we will provide the proof our main result, Theorem\,\ref{thm1.4}. Applying  the blow-up argument,  we are able to reduce the question of general atomic concentration problem  to that of one point concentration  in $\Bbb R^n$  treated in Section\,3, and applying Theorem 3.1 we conclude the proof. Finally, in Section\,5, using Corollary 1.2, we present the proof of Corollary \ref{cor1.3}.



\section{Local energy inequalities and the local pressure}
\label{sec:-2}
\setcounter{secnum}{\value{section} \setcounter{equation}{0}
\renewcommand{\theequation}{\mbox{\arabic{secnum}.\arabic{equation}}}}

In this section we  introduce the notion of  {\em local suitable weak solution} to the 
Euler equations similarly to the case of  the Navier-Stokes equations\cite{wol}.  As we shall prove below any solution satisfying Type I blow up  condition with  respect to the velocity  gradient is indeed local suitable weak solution before the possible blow-up time. 

\hspace{0.5cm}
Let us begin our discussion by recalling the definition of the local pressure  in a sub domain $ \Omega\subset \Bbb R^n $ with $ C^2$ boundary.  Here we distinguish between the two cases, firstly $ \Omega $ is bounded  and secondly, $ \Omega $ is an exterior domain. 

\vspace{0.2cm}
{\it 1.  Local pressure for $ \Omega $  bounded:} As in \cite{wol} we define the projection $ E^{ \ast}_\Omega : W^{-1,\, q}(\Omega ) \rightarrow W^{-1,\, q}(\Omega ) $ based on the unique solution of the Stokes equation as follows. Let $ f\in W^{-1,\, q}(\Omega )$ be given. Then we set $ E^{ \ast}_\Omega (f ):= \nabla p$, where $ p\in L^q_0(\Omega )$ stands for the unique pressure  from the 
unique weak solution  $ (w, p)\in W^{1,\, q}_{ 0, \sigma }(\Omega )\times L^q_0(\Omega )$ to the Stokes system 
\begin{align}
-\Delta w+ \nabla p &= f,\quad  \nabla \cdot w=0\quad  \text{in }\quad  \Omega ,\quad 
\label{2.1}
\\
w &=0 \quad  \text{ on}\quad \partial \Omega. 
\label{2.2b}
\end{align}
Here $ L^q_0(\Omega )$ stands for a subspace of all $ p\in L^q(\Omega )$ such that $ \intl_{\Omega } p dx=0$. 

(The existence and uniqueness in  bounded $ C^2$ domains  is due to Cattabriga\cite{catt},  while the case of bounded  $ C^1$   domains  were  treated in
\cite{gasiso}). Notice that $ \nabla p$ belongs to $ W^{-1,\, q}(\Omega )$ by
\[
 \langle \nabla p, \varphi \rangle := - \intl_{\Omega}  p \nabla \cdot \varphi dx, \quad  \varphi \in W^{1,\, q'}(\Omega ), {\color{blue} \quad 1/q+1/q'=1.}
\]
Obviously, from this definition it follows that $ E^{ \ast}(\nabla p)= \nabla p$ for every $ p\in L^q_0(\Omega )$, and thus it holds
\begin{equation}
(E^{ \ast}_\Omega) ^2 = E^{ \ast}_\Omega.  
\label{2.2a}
\end{equation} 
Observing the estimate 
\begin{equation}
\| \nabla w\|_{ L^q(\Omega )} + \| p\|_{ L^q(\Omega )} \le c \| f\|_{ W^{-1,\, q}(\Omega )},
\label{2.3}
\end{equation}
with a constant $ c>0$ depending only on $ q$ and the geometric property of $ \Omega $, we see that the operator $ E^{ \ast}_\Omega $ 
is bounded, satisfying
\begin{equation}
\| E^{ \ast}_\Omega (f)\|_{ W^{-1,\, q}(\Omega )} \le \| p\|_{ L^q(\Omega )} \le c \| f\|_{ W^{-1,\, q}(\Omega )} 
\label{2.4}
\end{equation}
with the same constant as in \eqref{2.3}.  

 In case $ f\in L^q(\Omega ) \hookrightarrow W^{-1,\, q}(\Omega )$, by virtue of the elliptic 
regularity of the Stokes system we find $E^{ \ast}_\Omega (f)= \nabla p \in L^q(\Omega ) $ together with the estimate   
\begin{equation}
\| \nabla p\|_{ L^q(\Omega )} \le c \| f\|_{ L^q(\Omega )},
\label{2.5}
\end{equation}
where $ c>0$ denotes a constant depending only on $ q$ and the geometric property of $ \Omega $. We also 
note  that in case $ \Omega $  equals to a ball $ B(x_0, r)$,  then  the constants  in both \eqref{2.4} and \eqref{2.5} depend neither on 
$ x_0$ nor on  $ r>0$.

\hspace{0.5cm}
In case $ 1 \le s \le +\infty$,  if  the vector valued function    $ f$ belongs to the Bochner space $ L^s(a,b; W^{-1,\, q}(\Omega ))$,  we may define $ E^{ \ast}_\Omega (f)$ pointwise 
\begin{equation}
E^{ \ast}_\Omega (f)(t) = E^{ \ast}_\Omega (f(t)),\quad  \text{ for a.e. $ t\in (a,b)$}. 
\label{2.5a}
\end{equation}
Clearly, \eqref{2.4} and \eqref{2.5}   imply that $ E^{ \ast}_\Omega $ is bounded on 
$ L^s(a,b; W^{-1,\, q}(\Omega ))$ and $ L^s(a,b; L^q(\Omega ))$ respectively.  For  $ f= \partial _t g$ in the sense of distributions  
we define
\[
  E^{ \ast}_\Omega (f) = \partial _t E^{ \ast}_\Omega (g)\quad  \text{ in the sense of distributions}.  
\]
 
 \vspace{0.2cm}
{\it 2.  Local pressure in case $ \Omega $  is an exterior domain:}  
Since $ \Omega $ is unbounded,  it will be more appropriate to replace the usual Sobolev space by the 
homogenous Sobolev space  $ D^{1,\, q}_0(\Omega )$,  which is defined as  the closure of $ C^{\infty}_{c}(\Omega )$ with respect to the norm 
\[
\| u\|_{ D^{ 1,2}(\Omega )} = \bigg(\intl_{\Omega } | \nabla u|^2 dx\bigg)^{ 1/2}. 
\]
Analogously, the subspace of all divergence free 
vector functions in  $ D^{1,\, 2}_{ 0}(\Omega; \R^n )$ will be denoted by $ D_{ 0, \sigma }^{ 1,2}(\Omega )$.  
In what follows by $ D^{ -1,2}(\Omega )$ we denote the dual of $ D^{ 1,2}_0(\Omega )$.  

\hspace{0.5cm}
As in the case of bounded domains for $ f\in D^{ -1,2}(\Omega )$ we define  $ E^{ \ast}_\Omega (f ) = \nabla p$, if 
$ (w, p)\in D^{ 1,2}_{ 0, \sigma }(\Omega; \R^n  )\times L^2(\Omega )$ denotes the unique weak solution to  Stokes problem 
\eqref{2.1}, \eqref{2.2b}.  Here the estimate \eqref{2.3} for $ q=2$ is still valid, which leads to the estimate 
\begin{equation}
\| E^{ \ast}_\Omega (f)\|_{ D^{-1,\, 2}(\Omega )} \le \| p\|_{ L^2(\Omega )} \le c \| f\|_{ D^{-1,\, 2}(\Omega )}. 
\label{2.6}
\end{equation}
This together with \eqref{2.2a} shows that $ E^{ \ast}_\Omega $ is a projection in $ D^{ -1,2}(\Omega )$ onto the closed subspace 
of all functionals $ \nabla p$ with $ p\in L^2(\Omega )$. In addition, if $ f\in L^{ q}(\Omega )$ for some $1< q < +\infty$,  then 
  $ \nabla p \in L^q(\Omega )$,  and there holds 
\begin{equation}
 \| \nabla p \|_{ L^q(\Omega )} \le c \| f\|_{ L^q(\Omega )}. 
\label{2.8}
\end{equation}
We also wish to remark that in case $ \Omega = B(x_0, r)^c = \R^n   \setminus B(x_0, r)$ in both \eqref{2.6} and \eqref{2.8}  
the  constants are independent of $ x_0$ and $ r>0$, which can be readily seen  by a standard  scaling argument.  
For vetor valued functions $ f\in L^s(a,b; D^{- 1,2}(\Omega ))$ we define $ E^{ \ast}_\Omega(f) $ and $ E^{ \ast}_\Omega(\partial _t f) $ as in the case of bounded domains.   

\hspace{0.5cm}
We are now in a position to introduce the notion of  local suitable weak solution to \eqref{1.1} in $ Q= \Omega \times (a,b)$.  

\begin{defin}
\label{def2.1}
A vector  function  $ v\in L^\infty(a,b; L^2(\Omega ))\cap L^3(Q)$ with $\nabla \cdot v=0 $ in the sense of distributions is said to be a 
{\it local suitable weak solution} to \eqref{1.1},  if the following two conditions are satisfied.

\hspace{0.5cm}
1. The function $ \widetilde{v}:= v + \nabla p_h: = v- E^{ \ast}_\Omega (v)$  solves 
\begin{equation}
\partial _t \widetilde{v} + (v\cdot \nabla) v = -\nabla p_0 \quad  \text{ in}\quad  Q
\label{2.9}
\end{equation}
 in the sense of  distributions, where $ \nabla p_h = -E^{ \ast}_{\Omega }(v)$ and $ \nabla p_0 = 
 -E^{ \ast}_{ \Omega }((v\cdot \nabla) v)$. 
 
 \hspace{0.5cm}
 2. For almost every $ a  \le  t < s < b$  and for all  $ \phi  \in C^{\infty}_{c}(\Omega )$ with $\phi \ge 0, $  the following local energy inequality holds true 
\begin{align}
&\intl_{\Omega } | \widetilde{v} (t)|^2 \phi   dx 
\cr
&\qquad \le   \intl_{\Omega } | \widetilde{v} (s)|^2 \phi dx  	
+  \intl_{t}^{s} \intl_{\Omega }  (|\widetilde{v}  |^2  v + 2 p_0 \widetilde{v} )\cdot \nabla \phi  dx d\tau 
 \cr
&\qquad \qquad + 
 \intl_{t}^{s} \intl_{\Omega }  v \otimes v: \nabla ^2 p_h  \phi  dx d\tau.
\label{2.10}
\end{align}
\end{defin}

\begin{rem}
\label{rem2.2} In \cite{shv}  the author has introduced the notion of  suitable weak solution under the assumption that the pressure $ p\in L^{ 3/2}(Q)$, and the local energy inequality holds true for almost every 
$ a  \le  t < s < b$  and for all  $ \phi  \in C^{\infty}_{c}(\Omega )$ with $\phi \ge 0, $
\begin{align}
&\intl_{\Omega } | v (t)|^2 \phi   dx 
 \le   \intl_{\Omega } | v (s)|^2 \phi dx  	
+  \intl_{t}^{s} \intl_{\Omega }  (|v |^2 +2p) v\cdot \nabla \phi  dx d\tau. 
 \label{2.12}
\end{align}
In fact, by the same the argument as in the proof of Lemma A.2 in \cite{ cw}, we see that any suitable weak solution satisfying \eqref{2.12} 
is also a local suitable weak solution in the sense of Definition\,\ref{def2.1}.  
\end{rem}

\hspace{0.5cm}
In the following lemma we show that any $ v$,  which satisfies the local energy inequality related to the  generalized energy inequality \eqref{2.10} for local suitable weak solutions,    admits  a weak measure valued trace in time.  

\begin{lem}
\label{lem2.3}
Let $ -\infty < a <b < +\infty$. Set $ Q= \Omega \times (a,b)$. 
Let $ v\in L^\infty(a,b; L^2(\Omega ))\cap L^3(Q)$, $ p\in L^{ 3/2}(Q)$,  $g\in L^{3}(Q)$,  and 
$f\in L^1 (a,b; L^2 (\Omega))$. 
Assume there exists a set $ J \subset (a,b]$ of Lebesgue measure zero such that the following local energy inequality holds 
true for all nonnegative $ \phi \in C^{\infty}_c(\Omega )$ and for all $ s,t \in (a,b]  \setminus J, t \le s$, 
 \begin{align}
& \intl_{\Omega } | v(t)|^2 \phi dx \le  \intl_{\Omega } | v(s)|^2 \phi dx  
\cr
& +   \intl_{t}^{s} \intl_{\Omega} (| v|^2 g + 2pv)\cdot \nabla \phi  dx d\tau + 
\intl_{t}^{s} \intl_{\Omega}   f\cdot v\phi  dx d\tau.
\label{2.13}
\end{align}
Then there exists a unique trace $ \sigma  \in L^\infty(a, b; \mathcal{M}^+(\Omega ))$ fulfilling the following properties:     
\begin{itemize}
  \item[(1)]  $ \sigma  (t) = | v(t)|^2 dx $ for a.e. $ t\in (a,b]$. 
  \item [(2)] The mapping $ t \mapsto \sigma  (t)$ is weakly-$ \ast$ left continuous, i.e. for every 
$ t\in (a,b]$ it holds 
\end{itemize}
\begin{equation}
\intl_{\Omega } \phi  d \sigma  (t) = \lim_{s \to t^-} \intl_{\Omega } \phi  d \sigma  (s)\quad  
\forall\,\phi \in C^{0}_{c}(\Omega ), \phi \ge 0.  
\label{2.14}
\end{equation} 
 \begin{itemize}
  \item[(3)] The following generalized local energy inequality holds for all $ a < t < s \le b$ and for all nonnegative 
  $ \phi \in C^{\infty}_c(\Omega )$
  
\end{itemize}
\begin{align}
& \intl_{\Omega }  \phi d\sigma  (t) \le  \intl_{\Omega }  \phi d\sigma  (s)  
\cr
& \qquad  +  \intl_{t}^{s} \intl_{\Omega} (| v|^2 g + 2pv)\cdot \nabla \phi  dx d\tau + 
\intl_{t}^{s} \intl_{\Omega} f\cdot v\phi dx d\tau.
\label{2.13a}
\end{align}\begin{itemize}
  \item[(4)]   The set $  \mathcal{M}_{ | v|^2}(b )$ contains only the measure $ \sigma _0 = \sigma (b)$. 
\end{itemize}

\end{lem}

{\bf Proof}: 
Let $ t\in (a,b]$. By $ \mathcal{M}(t)$ we define the set of all measures $ \sigma  \in \mathcal{M}(\Omega )$, obtained by a 
weak-$ \ast$ limit of the measures $ | v(\tau )|^2 dx $ as $ \tau \in [a,t)  \setminus J \rightarrow t$. In fact,  since $ | v(\tau )|^2 dx \in 
\mathcal{M}^+(\Omega )$ for all $ \tau \in (a,b)$,  we have $ \mathcal{M}(t) \subset \mathcal{M}^+(\Omega )$. 

 Let $ t_0\in (a,b]$,  and let 
$ \sigma  , \widetilde{\sigma  } \in \mathcal{M}(t_0)$ be two measures. Let $ \{ s_k\} $ be a sequence in $ (a, t_0)  \setminus J$ with $ s_k \rightarrow t_0$ as $ k \rightarrow +\infty$ such that 
\begin{equation}
| v(s_k )|^2 dx \rightarrow  \sigma  \quad  \text{ {\it weakly-$ \ast$ in } }\quad \mathcal{M}(\Omega ) \quad  \text{ {\it as}}\quad  k \rightarrow +\infty. 
\label{2.15}
\end{equation}   
From \eqref{2.13} with $ s=s_k$  we deduce,  after passing $ s_k \rightarrow t_0$,  that for all $ t\in (a,t_0)  \setminus J$ and for all nonnegative $ \phi \in C^{\infty}_c(\Omega )$
\begin{align}
& \intl_{\Omega } | v(t)| ^2\phi dx \le  \intl_{\Omega } \phi d \sigma    
\cr
& +   \intl_{t}^{t_0} \intl_{\Omega} (| v|^2g + 2pv)\cdot \nabla \phi  dx d\tau + 
\intl_{t}^{t_0} \intl_{\Omega}  f\cdot v \phi   dx d\tau.
\label{2.16}
\end{align}
Analogously, we take  a sequence  $ \{ \widetilde{s} _k\} $  in $ (a, t_0)  \setminus J$ with $ \widetilde{s} _k \rightarrow t_0$ as $ k \rightarrow +\infty$ such that 
\begin{equation}
| v(\widetilde{s} _k )|^2 dx \rightarrow  \widetilde{\sigma }  \quad  \text{ {\it weakly-$ \ast$ in } }\quad \mathcal{M}(\Omega ) \quad  \text{ {\it as}}\quad  k \rightarrow +\infty.  
\label{2.17}
\end{equation}   
Then from \eqref{2.16} with $ t=\widetilde{s} _k$  after passing $ \widetilde{s} _k \rightarrow t_0$ together with a standard mollifying  
argument we obtain 
\begin{equation}
\intl_{\Omega } \phi d \widetilde{\sigma }     \le \intl_{\Omega } \phi d \sigma \quad  \forall\,\phi \in C^{0}_c(\Omega ),\quad \phi  \ge 0.    
\label{2.18}
\end{equation}
Obviously, we may exchange $ \sigma  $ and $ \widetilde{\sigma  } $ in \eqref{2.18},  which yields the equality in \eqref{2.18}.  
Since  both $ \sigma  $ and $ \widetilde{\sigma  } $ are nonnegative  measures, we obtain  
 \[
\intl_{\Omega } \phi d \widetilde{\sigma }     = \intl_{\Omega } \phi d \sigma \quad  \forall\,\phi \in C^{0}_c(\Omega ).
\]
Thus $ \sigma  =\widetilde{\sigma  } $. This shows that for every $ t\in (a,b]$ there exists a unique nonnegative measure  
$ \sigma (t) \in \mathcal{M}^+(\Omega )$ such that 
 \begin{equation}
| v(\tau )|^2 dx \rightarrow  \sigma  (t)\quad  \text{ {\it weakly-$ \ast$ in } }\quad \mathcal{M}(\Omega ) \quad  
\text{ {\it as}}\quad  \tau  \rightarrow t^-. 
\label{2.19}
\end{equation}   
Furthermore, by the above definition of $ \sigma  (t)$ we get for all $ \phi \in C^0_c(\Omega )$ with  $ \max_{ \Omega}| \phi | \le 1$
\[
\intl_{\Omega } \phi  \sigma  (t) = \lim_{s_k \to t^-} \intl_{\Omega } | v(s_k)|^2 \phi dx 
\le \| v\|^2_{ L^\infty(a,b; L^2(\Omega ))},   
\]
which shows that $  \sigma \in L^\infty(a, b; \mathcal{M}^+(\Omega ))$.  

\hspace{0.5cm}
 In addition, from \eqref{2.19} we deduce that the following local energy inequality holds true for all $ s, t\in (a,b]$ with $ t \le s$ and for all nonnegative 
$ \phi \in C^{\infty}_c(\Omega )$ 
\begin{align}
& \intl_{\Omega } \phi \sigma  (t) \le  \intl_{\Omega } \phi d \sigma  (s)  
\cr
& +  \intl_{t}^{s} \intl_{\Omega} (| v|^2g + 2pv)\cdot \nabla \phi  dx d\tau + 
\intl_{t}^{s} \intl_{\Omega}  f\cdot v\phi   dx d\tau.
\label{2.20a}
\end{align}
By the same reasoning as the above it can be easily checked  that 
\begin{equation}
\sigma  (t) \rightarrow  \sigma  (s)\quad  \text{ {\it weakly-$ \ast$ in } }\quad \mathcal{M}(\Omega ) \quad  \text{ {\it as}}\quad 
t   \rightarrow s\quad  \text{ {\it in}}\quad (a, s). 
\label{2.21}
\end{equation}   
This implies that $\sigma :  t \mapsto\sigma  (t)$  is weakly-$ \ast$ left continuous, 
and therefore property $ (2)$ of the lemma  is fulfilled. To verify (1) of the lemma let $ t\in (a,b]$ be chosen so that 
\begin{equation}
\frac{1}{h}  \intl_{t-h}^{t}  | v(\tau )|^2 d\tau \rightarrow  | v(t)|^2\quad  \text{ {\it in}}\quad L^1(\Omega )\quad  \text{ {\it as}}\quad 
h \rightarrow 0^+. 
\label{2.22}
\end{equation}
As we have noted in Section\,1 due to Lebesgue's differentiation theorem for the Bochner integrable functions (see  e.g. \cite[Theorem 2, pp.134]{boc}) the property \eqref{2.22} holds true for a.e. $ t$.  It is also readily seen that from \eqref{2.22} 
 we get for all $ \phi \in  C^{0}_c(\Omega )$
 \begin{equation}
 \frac{1}{h}  \intl_{t-h}^{t}  \intl_{\Omega}  | v(\tau )|^2 \phi dx d\tau \rightarrow   \intl_{\Omega}  | v(t )|^2 \phi dx\quad  \text{ as}\quad 
h \rightarrow 0^+. 
 \label{2.23}
 \end{equation} 
 We fix $ \phi \in C^{0}_c(\Omega )$. 
  Let $ \{ h_k\}$ be a sequence in $ (0, t-a)$ which 
 converges to zero as $ k \rightarrow +\infty$.  By the mean value theorem  for the integrals for every $ k\in \N$ we may choose $ t_k \in (t- h_k, t)  \setminus J$ such that 
\[
  \intl_{\Omega}  | v(t_k)|^2 \phi dx=  \frac{1}{h_k}  \intl_{t-h_k}^{t}  \intl_{\Omega}  | v(\tau )|^2 \phi dx d\tau.
\]
This together with \eqref{2.23} and the weakly-$ \ast $ left side continuity of $ \sigma $ yields 
\[
  \intl_{\Omega} \phi d\sigma  (t) = \lim_{k \to \infty}  \intl_{\Omega}  | v(t_k)|^2 \phi dx = \intl_{\Omega}  | v(t )|^2 \phi dx,
\] 
and therefore (1) of the lemma  is satisfied.  

\hspace{0.5cm}
Finally, the generalized local energy inequality \eqref{2.13a}   follows immediately from \eqref{2.13} together 
with \eqref{2.19}, while (4) of the lemma  immediately follows from the proof of (1). In fact, we already have proved that 
$ \sigma (t)= | v(t)|^2 dx$ for every  Lebesgue point of $ | v(\cdot )|^2$, which immediately gives (4),  since for every 
$| v(t)|^2 dx \rightarrow \sigma (0) $ as $ t \rightarrow b^-$ for $ t$ in the Lebesgue set of $ | v(\cdot )|^2$ in $ (a,b)$. 

\hfill \Beweisende

As an  important consequence of Lemma\,\ref{lem2.3} we are able to study the concentration for the  local suitable 
weak solutions to the Euler equation. In fact we have the following.

\begin{rem}
\label{rem2.4}   1. If  $ v\in L^\infty(a, b, L^2(\Omega ))\cap L^3(Q)$ is a local suitable weak solution  to the Euler equations in 
$Q= \Omega \times (a,b)$,  then $ \widetilde{v} = v +\nabla p_h$  is a distributional solution to 
\begin{equation}
\partial _t \widetilde{v} + (v\cdot \nabla) v = - \nabla p_0\quad  \text{ in}\quad  Q,
\label{eulerLP}
\end{equation}
where 
\[
\nabla p_{ h} = - E^{ \ast}_{ \Omega }(v),\quad  \nabla p_0= -E^{ \ast}_\Omega ((v\cdot \nabla) v). 
\]
Note that
 \[
 (v\cdot \nabla) v = (v \cdot \nabla) \widetilde{v} - v\cdot \nabla^2 p_h.
\]
Since $ \widetilde{v}$ fulfills \eqref{2.10},    the local energy inequality \eqref{2.13} holds for $ \widetilde{v} $ 
 in place of $ v$  for a.e. $ a< t <s<b$ with 
$$ g= v,\quad  f=v\cdot \nabla^2 p_h. $$
 According to Lemma\,\ref{lem2.3} 
there exists a unique $ \widetilde{\sigma }\in L^\infty(a, b;  \mathcal{M}^+ (\Omega ))$   such that $ (1)$-$ (4)$ of the lemma are fulfilled. 
In particular, we see that $ M_{ | \widetilde{v} |^2}(b)= \{ \widetilde{\sigma}(b)  \}$, and there holds 
\begin{equation}
\begin{cases}
| \widetilde{v}  (t)|^2 dx \rightarrow \widetilde{\sigma }(b)  \quad  \text{{\it weakly-$ \ast$ in }}\quad \mathcal{M}(\Omega )  \quad  
  \\[0.3cm]
\text{{\it as}}\quad  t \rightarrow b^-\quad  \text{ {\it for $ t$ chosen from  the Lebesgue set of \,\, $ | \widetilde{v}(\cdot ) |^2$}}.
\end{cases}
\label{2.27}
\end{equation}
While  the set $ \mathcal{M}_{ | \widetilde{v} |^2}(b)$ contains only one unique measure, it is not  true  in general for 
$ \mathcal{M}_{ |v|^2}(b)$. The reason is that $ v$ may not satisfy the local energy inequality.  However, as we shall  show below 
by Lemma\,\ref{lem2.4} the concentration set of measures in $ \mathcal{M}_{ | v|^2}(b)$ coincides with the concentration set of $ \widetilde{\sigma}  (b)$, 
which is the unique measure in $ \mathcal{M}_{ | \widetilde{v} |^2}(b)$. 
 
\hspace{0.5cm}
2. In case that $ v$ is a solution to the Euler equations \eqref{1.1} satisfying  Type I blow-up condition with 
respect to the velocity gradient,  then $ v$ is a local suitable weak solution in the sense of Definition\,\ref{def2.1}. In other words, 
$ \widetilde{v} = v+ \nabla p_{ h}$ satisfies the energy inequality \eqref{2.10} for all $ s,t\in [-1,0)$, $ t \le s$.  
As we mentioned above, thanks  to Lemma\,\ref{lem2.3} there exists a unique measure valued trace 
$ \widetilde{\sigma  } \in \mathcal{M} (\Omega )$. Since every $ t\in (0,1)$  is a Lebesgue point of $ | v|^2$ it follows 
$ \widetilde{\sigma} (s)= | v(s)|^2 dx $  for all $ s\in (-1,0)$, and there holds  
\begin{equation}
 | \widetilde{v} (s)|^2 dx \rightarrow \widetilde{\sigma  }(b)   \quad  \text{{\it weakly-$ \ast$ in}}\quad \mathcal{M}(\Omega )  
 \quad \text{{\it as}}\quad  s \rightarrow 0^-.  
\label{2.24}
\end{equation}
 
\end{rem}

\begin{lem}
\label{lem2.4}
Let $ v\in L^\infty(a,b; L^2(\Omega ))\cap L^3(Q)$ be a local suitable weak solution to \eqref{1.1}.  Let $ \nabla p_{ h}$ and $ \nabla p_0$ 
denote the corresponding local pressure (cf. Definition\,\ref{def2.1}), and define    $ \widetilde{v}:= v+ \nabla p_h $. Then 
each measure $ \sigma _0 \in \mathcal{M}_{ | v|^2}(b)$ has no atoms if and only if $ \widetilde{\sigma }(b) $ has no atoms.  
\end{lem}

{\bf Proof}: Let $ \sigma _0 \in \mathcal{M}_{ | v|^2}(b)$,  then there exists a sequence $\{ s_k\} $ 
in the set of  Lebesgue points 
of $| v(\cdot )|^2 $  such that $ s_k \rightarrow b^-$ and 
\begin{equation}
 | v (s_k)|^2 dx \rightarrow \sigma _0   \quad  \text{{\it weakly-$ \ast$ in}}\quad \mathcal{M}(\Omega )   \quad  
  \text{{\it as}}\quad  k \rightarrow +\infty.  
\label{2.30}
\end{equation}
Since  $ \{ v(s_k)\}$ is bounded in $ L^2(\Omega )$, thanks to  the reflexivity, eventually passing to a subsequence,  we may assume that there exists  
$ v_0\in L^2(\Omega )$ such that  
\[
 v(s_k) \rightarrow v_0  \quad  \text{{\it weakly in}}\quad  L^2(\Omega )\quad  \text{{\it as}}\quad  k \rightarrow +\infty.  
\] 
 By the boundedness of the operator $ E^{ \ast}_{ \Omega } $ in $ L^2(\Omega )$  we deduce that  
\[
 \nabla p_h(s_k) \rightarrow \nabla p_{ h,0}  \quad  \text{{\it weakly in}}\quad  L^2(\Omega )\quad  \text{{\it as}}\quad  k \rightarrow +\infty,
\] 
 where $ \nabla p_{ h, 0}= - E^{ \ast}_\Omega (v_0)$.   By virtue of Lemma\,\ref{lemA.2} we find  that  for every $ \Omega ' \Subset \Omega $ the 
 convergence  $\nabla p_h(s_k) \rightarrow \nabla p_{ h,0}$ as $ k \rightarrow +\infty$ is in fact uniform on $ \Omega '$. This 
 together with  the weak convergence  of $ v(s_k) \rightarrow v_0$ in $ L^2(\Omega )$ implies 
 that     
 \begin{equation}
 \begin{cases}
\Big(2 v(s_k)\cdot \nabla p_{ h}(s_k)  - | \nabla p_h (s_k)|^2\Big) dx \rightarrow
\Big(2 v_0\cdot \nabla p_{ h,0} - | \nabla p_{ h,0}|^2\Big) dx 
 \\[0.3cm]
  \text{{\it weakly-$ \ast$ in }}\quad   \mathcal{M}(\Omega )\quad  
  \text{{\it as}}\quad  k \rightarrow +\infty.  
 \end{cases}
 \label{2.31}
\end{equation}
Furthermore, verifying that $ \widetilde{v}\in C_w([a,b]; L^2(\Omega ))$ and recalling the  weakly-$ \ast$ left side 
continuity of $ \widetilde{\sigma } $,  we have deduce that 
\begin{equation}
| \widetilde{v} (t)|^2 dx \le \widetilde{\sigma}  (t) \quad  \forall\,t\in (a,b]. 
\label{2.32}
\end{equation}
Combining \eqref{2.30} and \eqref{2.31}, noting $ | v(s_k)|^2 = | \widetilde{v} (s_k)|^2 - 
2v(s_k)\cdot \nabla p_{ h}(s_k)+ | \nabla p_{ h}(s_k)|^2$, and employing \eqref{2.32},    we infer that  
\begin{align}
0 \le \sigma_0 &\le \text{$ w^{ \ast}$\!-\!}\lim_{k \to \infty}\Big( \widetilde{\sigma }(t_k)  + (- 2 v_0 \cdot \nabla p_{ h,0} +  | \nabla p_{ h,0}|^2) dx\Big) 
\cr
&= \widetilde{\sigma}  (b) + \Big( - 2 v_0 \cdot \nabla p_{ h,0} +  | \nabla p_{ h,0}|^2)\Big) dx. 
\label{2.33}
\end{align} 
This immediately shows that if $ \widetilde{\sigma} (b)$ has no atoms, the same also holds true for  $\sigma _0 $. 

\hspace{0.5cm}
In order to prove the opposite direction  we assume that each measure in  $ \mathcal{M} _{ | v|^2}(b)$ has no atoms. 
Let us choose a sequence $ \{ s_k\} $ in $ (a,b)$ such that $ s_k \rightarrow b$ 
as $ k \rightarrow +\infty$ with the property that each $ s_k$ is simultaneously belong to the Lebesgue set of 
 $ | v(\cdot )|^2$ and $ | \widetilde{v}(\cdot ) |^2$.  Eventually passing to a subsequence,  we may assume there exist a measure 
$ \sigma _0 \in \mathcal{M}^+(\Omega )$ and $ v_0\in L^2(\Omega )$   having the following convergence properties 
\begin{align}
&  | v (s_k)|^2 dx \rightarrow \sigma _0   \quad  \text{{\it weakly-$ \ast$ in}}\quad \mathcal{M}(\Omega ),   
\label{2.34}
\\
&v (s_k) \rightarrow v_0  \quad  \text{{\it weakly in}}\quad  L^2(\Omega )\quad  \text{{\it as}}\quad  k \rightarrow +\infty.  
\label{2.35}
\end{align} 
Thanks to the property (4) of Lemma\,\ref{lem2.3} it holds  
\begin{equation}
 | \widetilde{v}  (s_k)|^2 dx \rightarrow \widetilde{\sigma}  (b)   \quad  \text{{\it weakly-$ \ast$ in}}\quad \mathcal{M}(\Omega )
\quad  \text{{\it as}}\quad  k \rightarrow +\infty.  
\label{2.36}
\end{equation}
Arguing as in the first part of the proof,  from \eqref{2.34} and \eqref{2.35} we get the  property \eqref{2.31}.  
Finally, observing \eqref{2.36}, we conclude that 
\begin{align*}
\sigma_0 &= \text{$ w^{ \ast}$\!-\!}\lim_{k \to \infty}  |v(s_k) |^2 dx
\\
&=\text{$ w^{ \ast}$\!-\!}\lim_{k \to \infty} \Big(| \widetilde{v}(s_k) |^2 - 2 v(s_k)\cdot \nabla p_{ h}(s_k)  + | \nabla p_h (s_k)|^2\Big) dx 
\\
&= \widetilde{\sigma}  (b)  + \Big(- 2 v_0 \cdot \nabla p_{ h,0} +  | \nabla p_{ h,0}|^2\Big) dx. 
\end{align*} 
Since $ \sigma _0$ has no atoms,  the above identity  shows  that $ \widetilde{\sigma }(b) $ also has no atoms.  \hfill \Beweisende

\section{Removing one point energy  concentration in $\Bbb R^n$}
\label{sec:-3}
\setcounter{secnum}{\value{section} \setcounter{equation}{0}
\renewcommand{\theequation}{\mbox{\arabic{secnum}.\arabic{equation}}}}

In this section we restrict ourself to the case $ \Omega = \R^n $. In this case,  since  any solution which  satisfies Type I condition with 
respect to the velocity gradient enjoys the local energy inequality, the pressure  satisfies $ p\in L^{ 3/2}( \R^n \times (-1,0))$ due to the Calder\'on-Zygmund inequality, and 
 thanks to Lemma\,\ref{lem2.3} there exists a unique measure $ \sigma  \in \mathcal{M}(\R^n )$ such that 
\begin{equation}
 |v(t) |^2 dx \rightarrow \sigma    \quad  \text{{\it weakly-$ \ast$} in}\quad  \mathcal{M}(\R^n )\quad  \text{{\it as}}\quad  t \rightarrow 0^-.  
\label{conc2}
\end{equation}
    
Our aim  is the proof of  Theorem\,\ref{thm1.4} for the special case that $ \sigma  $ in \eqref{conc2} equals to
the Dirac measure $ E _0 \delta _0$ for some constant $ 0 \le E_0 < +\infty$. Namely we shall prove the following:

\begin{thm}
\label{thm3.1}
Let $ v\in L^2(-1, 0; L^2_\sigma (\R^n ))$ be a solution to the Euler equations \eqref{1.1}.  In addition, we assume 
that $ v$ satisfies the Type I blow up condition \eqref{rate} (cf. Theorem\,\ref{thm1.4}) and \eqref{conc2} with   
$\sigma  = E_0   \delta _0$ for some $ 0 \le E_0 < +\infty$. Then $ v \equiv 0$.   
\end{thm}
  
 \begin{rem}
 In the proof of Theorem\,\ref{thm3.1}   
 we make significant use of several decay properties of the solution to the Euler equations with respect to the  space and  time variables as we approach the blow-up time. The decay estimate  is  actually obtained under following  more general condition than \eqref{rate}
\begin{equation}
\exists\,\mu \in \Big[\frac{n}{n+2},1\Big):\quad   \sup_{ t\in (-1,0)} (-t)^{ \frac{n+2}{n}\mu } \| \nabla v(t)\|_{ L^\infty} <+\infty. 
\label{rate1}
\end{equation}
 \end{rem}

We divide the proof of Theorem\,\ref{thm3.1} into  four steps, each step being a subsection below.

\subsection{Proof of $  E_0=E$}

The aim of this section is to show that  $ E_0 = E = \| v(t)\|^2_{ L^2}$ , $ -1 \le t <0$,  in Theorem\,\ref{thm3.1} under the  condition \eqref{conc2},
 in other words,  the energy cannot escape into infinity at the blow-up time.
We begin with the following observation.

\begin{lem}
\label{lem3.3}
Let $ v\in L^\infty(-1, 0; L^2(\R^n ))$,  which satisfies  \eqref{rate1}.  Then, 
 it holds 
\begin{equation}
\sup_{ t\in (-1, 0)} (-t)^{ \a } \| v(t)\|_{ L^\infty} \le  c ( E^{ \frac{1}{2}} + \sup_{ t\in (-1, 0)}(-t)^{ \frac{(n+2)\mu}{n}}\| \nabla v\|_{ L^\infty})=: C_1<+\infty,
\label{1.3d}
\end{equation} 
 where $ E= \| v(-1)\|^2_{ L^2}$. 
\end{lem}

{\bf Proof}:  This is immediate of the Gagliardo-Nirenberg inequality and the energy conservation $E=E(t)$ for $t\in (-1, 0)$(see Remark \ref{rem1.4}),
\begin{align*}
(-t)^{ \a } \| v(t)\|_{ L^\infty}&
\leq c(-t)^{ \a } \|v(t)\|_{L^2}^{\frac{2}{n+2}} \|\nabla v(t)\|_{L^\infty}^\frac{n}{n+2}\leq c E^{\frac{1}{n+2}} \left\{(-t)^{  \frac{(n+2)\mu}{n}}\|\nabla v(t)\|_{L^\infty}\right\}^\frac{n}{n+2}\\
&\le   c   (E^{\frac12} + (-t)^{ \frac{(n+2)\mu}{n}}\| \nabla v\|_{ L^\infty}).
\end{align*}
  \hfill \Beweisende 
  
From Lemma\,\ref{lem3.3} along with $ v\in L^\infty(-1,0; L^2(\R^n ))$  we immediately get 
\begin{equation}
v\in L^1(-1,0; L^\infty(\R^n ))\cap L^3(-1,0; L^3(\R^n )). 
\label{Linfty}
\end{equation}

We have the following 

\begin{lem}
\label{lem3.4}
Let $ v\in L^\infty(-1,0; L^2_{ \sigma }(\R^n ))\cap L^\infty_{ loc}([-1,0), W^{1,\, \infty}(\R^n ))$ be a solution  to \eqref{1.1}
satisfying \eqref{rate1} and  \eqref{conc2} with $ \sigma  = E_0   \delta _0$. Then
it holds $ E_0=E$. 
\end{lem}

{\bf Proof}: Given  $ 0< R< +\infty$,  we denote  by  $ \eta _R\in C^{\infty}_c(\R)$ a cut off function such that 
$ 0 \le \eta_R \le  1$ in $ \R$, $ \eta_R = 1$ on $ (1, R)$,  $ \eta_R =0$ in $ (-\infty, 0]\cap (2R, +\infty)$ and 
$ | \eta'_R| \le \frac{2}{R}$ in $ (R, 2R)$.    We multiply  the Euler equations by $ - v \eta_R(| x|^2)$,  integrate the result 
over $ \R^n \times (t, \tau ), -1 \le t < \tau  < 0$, and apply integration by parts. This gives 
\begin{align*}
&\frac{1}{2} \intl_{ \R^n } | v(t)|^2 \eta_R (| x|^2) dx 
\\
&\qquad = \frac{1}{2} \intl_{ \R^n } | v(\tau )|^2 \eta_R (| x|^2) dx + \frac{1}{2}  \intl_{t}^{\tau }  
\intl_{ \R^n } | v(s)|^2 v(s )\cdot \nabla  \eta_R (| x|^2) dxds 
\\
&\qquad \qquad  +  \intl_{t}^{\tau }  \intl_{ \R^n } p(s) v(s)\cdot \nabla  \eta_R (| x|^2) dxds.
\end{align*}
Observing \eqref{conc2},   we see that  $  \intl_{ \R^n } | v(\tau )|^2 \eta_R (| x|^2) dx  \rightarrow 0$
as $ \tau \rightarrow 0$. In view of \eqref{Linfty} together with the Calder\'on-Zygmund estimate,  having   $ v\in L^3(\R^n \times (-1,0))$ and $ p\in L^{ 3/2}(\R^n \times (-1,0))$, we obtain from the above identity after letting $ \tau  \rightarrow 0$
\begin{align*}
&\frac{1}{2} \intl_{ \R^n } | v(t)|^2 \eta_R (| x|^2) dx 
\\
&\qquad =\frac{1}{2}  \intl_{t}^{0}  
\intl_{ \R^n } | v(s)|^2 v(s )\cdot \nabla  \eta_R (| x|^2) dxds 
 +  \intl_{t}^{0}  \intl_{ \R^n } p(s) v(s)\cdot \nabla  \eta_R (| x|^2) dxds.
\end{align*}
We are now in a position to pass $ R \rightarrow +\infty$  in the above to get 
\begin{align*}
&\intl_{ \R^n } | v(t)|^2 \eta(| x|^2) dx 
\\
&\qquad = \intl_{t}^{0}  
\intl_{ \R^n } | v(s)|^2 v(s )\cdot \nabla  \eta(| x|^2) dxds 
 +  2\intl_{t}^{0}  \intl_{ \R^n } p(s) v(s)\cdot \nabla  \eta(| x|^2) dxds,
\end{align*}
where $ \eta \in C^{\infty}(\R)$ stands for the corresponding cut off function such that $ \eta \equiv 1$ on $ (1, +\infty)$. 
Noting that $ 1- \eta (| x|^2)\in C^{\infty}_{c}(\R^n )$ and once more appealing to \eqref{conc2},  
from the above identity we deduce 
\begin{align*}
E  &= \intl_{ \R^n } | v(t)|^2 \eta(| x|^2) dx 
+ \intl_{ \R^n } | v(t)|^2 \left(1-\eta(| x|^2) \right) dx  
\\
&= \intl_{ \R^n } | v(t)|^2 \left(1-\eta(| x|^2) \right) dx   +  \intl_{t}^{0}  
\intl_{ \R^n } | v(s)|^2 v(s )\cdot \nabla  \eta(| x|^2) dxds 
\\
&\qquad  +  2\intl_{t}^{0}  \intl_{ \R^n } p(s) v(s)\cdot \nabla  \eta(| x|^2) dxds,
\\
&\rightarrow E_0  \quad  \text{ as}\quad  t \rightarrow 0.  
\end{align*}
 
Whence, the claim.  \hfill \Beweisende

\subsection{Decay estimates for energy concentrating solutions}

In this subsection our aim is to prove the  space-time decay for solutions to the Euler equations  satisfying the blow-up rate \eqref{rate1} and the energy concentration at 
$ (0,0)$.

\begin{lem}\label{spacedecay}
Let $ v\in L^2(-\infty, 0; L^2_\sigma (\R^n ))\cap L^\infty_{ loc}([-1,0), W^{1,\, \infty}(\R^n ))$ be a  solution to the Euler equations satisfying \eqref{rate1} and \eqref{conc2} with $ \sigma  = E  \delta _0$. 
Then for every $  0 < \beta < n+2$ there exists a constant $ c$ depending on $ C_1, \mu$ and $ \beta $ such that  
for every $ t\in [-1, 0)$ it holds 
\begin{equation}
\intl_{\R^n } | v(t)|^2  | x|^{ \beta } dx \le c (-t)^{ \beta(1-\a ) }. 
\label{3.1a}
\end{equation}
\end{lem}

{\bf Proof}: 
Given $R>2$, let $ \eta_R \in C^{\infty}_{\rm c}(\R)$  
denote a cut off function such that 
$ 0 \le  \eta_R \le 1$ in $ \R$, $ \eta_R \equiv 0$ in $ (2R, +\infty)$, 
$ \eta_R \equiv 1 $ in $ (-\infty, R)$, and  $ | \eta_R '| \le \frac{2}{R} $ in $\R$.  
Observing  \eqref{conc2} with $ \sigma = E  \delta _0$, we get for all $ 0< \beta < +\infty$
$$ 
\lim_{t\to 0-} \int_{\Bbb \R^n } |v(t)|^2  \eta_R (|x| )^2| x|^\beta dx=  E  \Big\langle  \delta _0,  \eta_R (|x| )^2| x|^\beta \Big\rangle=0. 
$$
 We multiply (1.1) by 
$ - v \eta _R(| x|) ^2| x|^\beta $, $ 1 \le  \beta < +\infty$,   integrate over $ \R^n \times (t,0)$,  $ -1 < t< 0$, and apply integration by parts. 
This together with $ \nabla \Big( \eta _R(| x|)^2 | x|^\beta\Big)  =  \beta  x | x|^{ \beta -2}\eta _R(| x|) +2x \eta _R(| x|) \eta '_R(| x|)| x|^{ \beta -1}  $  gives

\begin{align}
&\frac{1}{2} \intl_{\R^n } | v(t)|^2 \eta_R(| x|)^2| x|^\beta   dx  \cr
&\qquad=  - \frac{\beta }{2}\intl_{t}^{0} \intl_{\R^n } v(s)\cdot x | x|^{ \beta -2} | v(s)|^2  \eta _R(| x|)^2 dxds
\cr
 &\qquad \qquad-  
\intl_{t}^{0} \intl_{\R^n } v(s)\cdot x\,\eta _R(| x|)\eta' _R(| x|) | x|^{ \beta-1 }| v(s)|^2  dxds 
\cr
& \qquad \qquad- \beta \intl_{t}^{0} \intl_{\R^n } p(s) v(s)\cdot x | x|^{ \beta -2}  \eta _R(| x|)^2 dxds 
\cr
&\qquad \qquad- 
2 \intl_{t}^{0} \intl_{\R^n } p(s) v(s)\cdot x\, \eta _R(| x|)\eta' _R(| x|) | x|^{ \beta -1}   dxds 
\cr
& = I+II+III+IV.
\label{3.1b}
\end{align}

\hspace{0.5cm}
In what follows, we will make an extensive use of the following estimate  
\begin{equation}
\intl_{\R^n } | p(s)|^2 | x|^{ \gamma } dx \le c\intl_{\R^n } |  v(s)|^4 | x|^\gamma  dx  \le 
 c(-s)^{ -2\a} \intl_{\R^n } |  v(s)|^2 | x|^\gamma  dx
\label{3.1gg}
\end{equation} 
which holds true for all $ 0 \le \gamma < n$.
Indeed, in case $ \gamma =0$ the estimate   \eqref{3.1gg} is an immediate consequence of the well-known Calder\'on-Zygmund inequality together with \eqref{1.3d}. 
For $ 0< \gamma < n$, noting that $ | x|^\gamma  $  belongs to the class $ A_2$,  the estimate \eqref{3.1gg} follows by the aid of the weighted Calder\'on-Zygmund  inequality  \cite[Corollary, p.205]{ste}  along with  \eqref{1.3d}.  

\hspace{0.5cm}
We divide the proof of \eqref{3.1a} into five steps:  

\hspace{0.5cm}
{\it 1. We consider the case $ \beta =1$}: 
Noting that $ \eta '_R(| x|) | x| \le 4$ and observing \eqref{1.3d},  we immediately get 
\[
I+ II \le 5c \intl_{t}^{0} (-s)^{ -\a } \intl_{\R^n }   | v(s)|^2  dxds \le c   E  (-t)^{1-\a }.
\]

For $ s\in (-1, 0)$ observing that $  \Delta  p (s) = -\nabla \cdot \nabla \cdot  (v(s) \otimes v (s))$, the estimate  \eqref{3.1gg} for $ \gamma =0$ 
gives 
\[
\| p(s)\|_{ L^2}  \le c(-s)^{ -\a }  \| v(s)\|_{ L^2}\le c E^{ \frac{1}{2}} (-s)^{ -\a } .  
\] 
The above inequality along with Cauchy-Schwarz's inequality yields

\begin{align*}
III+ IV &\le   9\intl_{t}^{0}   \intl_{\R^n } | p(s)| | v(s)|dxds \le c  \intl_{t}^{0} \| p(s)\|_{ L^2} \| v(s)\|_{ L^2}   ds
\\
&\le c E^{ \frac{1}{2}}  \intl_{t}^{0} (-s)^{ -\a }  ds = c  E^{ \frac{1}{2}} (-t)^{1-\a }.
\end{align*}
Hence, from \eqref{3.1b} it follows that  
\[
 \intl_{\Bbb R^2} | v(t)|^2 \eta _R(| x|)^2 | x| dx  \le c  E^{ \frac{1}{2}} (-t)^{1-\a }. 
\]

After passing $ R \rightarrow +\infty$  in the above inequality, we get the estimate   \eqref{3.1a} for $ \beta =1$.  

\vspace{0.2cm}
{\it 2. We consider the case $ \beta =2$: }
Noting, that $ \| v(s)\|_{ L^\infty} \le  c(-s)^{ -\a } $ and 
$ | \eta' _R| | x| \le 4$ together with \eqref{3.1a} for $ \beta =1$, we easily find 
\[
I+ II \le 4 \intl_{t}^{0}(-s)^{ -\a }\intl_{\R^n }   | v(s)|^2  | x| dxds \le c (-t)^{2-2\a}.
\]

Applying \eqref{3.1gg} for $  \gamma =1$ and making use of  \eqref{3.1a} for  $\beta=1  $, we get 
\[
\intl_{\R^n } | p(s)|^2 | x| dx \le c(-s)^{ -2\mu }  \intl_{\R^n } |  v(s)|^2 | x|  dx
  \le  c (-s)^{ 1-3\a}. 
\]
Integrating the above inequality over $ (t,0)$,  and using once more  \eqref{3.1a} for $ \beta =1$, we obtain 
\begin{align*}
III + IV &\le c \int_t ^0 \left(\int_{\Bbb \R^n } |p(s)|^2 |x|dx \right)^{\frac12} \left(\int_{\Bbb \R^n } |v(s)|^2 |x|dx \right)^{\frac12} ds\\
 &\le c\intl_{t}^{0}(-s)^{ \frac12-\frac{3\a}{2}} \bigg( \intl_{\R^n } | v(s)|^2 | x| dx\bigg)^{ 1/2}  ds\\
&\le c\intl_{t}^{0} (-s)^{ 1-2\mu }ds = c (-t)^{2-2\mu }. 
\end{align*}
Inserting the estimates of $ I, II, III$, and $ IV$ into \eqref{3.1b}, and passing $R\to +\infty$,  we obtain 
\begin{equation}
\intl_{\R^n } | v(t)|^2  | x|^2 dx \le c (-t)^{2(1-\a )}. 
\label{3.1c}
\end{equation}
 
 {\it 3.} Iterating  the above argument for $ \beta =3, \ldots, n$, by using \eqref{3.1gg} for $ 0<\gamma <n$, we find 
\begin{equation}
\intl_{\R^n } | v(t)|^2  | x|^n dx \le c (-t)^{ n(1-\a ) } . 
\label{3.1d}
\end{equation}

{\it 4. Next, we we consider the case $ \beta =n+1$}. Arguing as above, in this case  we estimate  
\[
I+II \le c (-t)^{ ( n+1)(1-\a )}. 
\] 
For the estimation of $ III$ and $ IV$ we make use of \eqref{3.1gg} for $ \gamma =n-1$, Cauchy-Schwarz' inequality and Young's inequality  to get 
\begin{align*}
III+ IV  &\le  c  \intl_{t}^{0} \intl_{ \R^n }   |  p(s)| | v(s)|\eta _R(| x|) | x|^n dxds
\\
&\le  \intl_{t}^{0}      \bigg( \intl_{\R^n } | p|^2 | x|^{ n-1 }  dx\bigg)^{ 1/2}  \bigg(\intl_{\R^n } | v(s)|^2  \eta _R(| x|)^{ 2}| x|^{ n+1} dx\bigg)^{ 1/2}ds
\\
&\le \intl_{t}^{0}  (-s)^{-\mu+ \frac{(n-1)(1-\mu )}{2} } \bigg(\intl_{\R^n } | v(s)|^2  \eta _R(| x|)^2| x|^{ n+1} dx\bigg)^{ 1/2} ds\\
&\le  c(-t)^{1-\mu+\frac{(n-1)(1-\mu )}{2} }  \bigg(\esssup_{s\in (t, 0)}\intl_{\R^n } | v(s)|^2  \eta _R(| x|)^2| x|^{ n+1} dx\bigg)^{ 1/2} 
\\
&\le c(-t)^{ (n+1) (1-\mu)} + \frac{1}{4}\esssup_{s\in (t, 0)} \intl_{\R^n } | v(s)|^2  \eta _R(| x|)^2 | x|^{ n+1} dx.
\end{align*}
Inserting the above estimates of $ I, II, III$, and $ IV$ into \eqref{3.1b},  the following inequality holds for all 
$-1 <  t   <0$  
\begin{align*}
& \intl_{\R^n } | v(t)|^2  \eta _R(| x|)^2 | x|^{ n+1} dx 
\\
&\le c(-t)^{ (n+1) (1-\mu)} + \frac{1}{2}\esssup_{s\in (t, 0)} \intl_{\R^n } | v(s)|^2  \eta _R(| x|)^2 | x|^{ n+1} dx.
\end{align*}
Let $ -1 < \tau < 0$. Taking supremum over $ t\in (\tau , 0) $ in both sides of the above inequality and noting that 
function on the right-hand side attains the maximum at $ t=\tau $, we get  
\begin{align*}
&\esssup_{t\in (\tau , 0)} \intl_{\R^n } | v(t)|^2  \eta _R(| x|)^2 | x|^{ n+1} dx 
\\
&\le c(-\tau )^{ (n+1) (1-\mu)} + \frac{1}{2}\esssup_{s\in (\tau , 0)} \intl_{\R^n } | v(s)|^2  \eta _R(| x|)^2 | x|^{ n+1} dx.
\end{align*}
Accordingly,  for all $ -1< t< 0$ it holds 
 \begin{equation}
\intl_{\R^n } | v(t)|^2  | x|^{ n+1} dx \le c (-t)^{ (n+1)(1-\mu )}. 
\label{3.1f}
\end{equation}

{\it 5. We now consider the case $ n+1 < \beta < n+2$.} Using H\"olders inequality  together with $ \| v(s)\|_{ L^2}= E^{ \frac{1}{2}}$ and \eqref{3.1f}, we deduce that   for all $ 0 < \gamma \le n+1$ and $ -1 \le s < 0$ it holds 
\begin{equation}
\intl_{\R^n } | v(s)|^2  | x|^{ \gamma } dx \le c (-s)^{\gamma(1-\mu )  }. 
\label{3.1g}
\end{equation}
Applying  the estimate \eqref{3.1gg} for  $ 0< \gamma< n $, and using \eqref{3.1g},  we get 
\begin{equation}
\intl_{ \R^n } | p(s)|^2 | x|^{ \gamma } dx   \le  c
(-s)^{\gamma (1-\mu )  -2\mu }.
\label{3.1h}
\end{equation}
We now easily estimate $I+ II $ by using \eqref{3.1g} with $ \gamma = \beta -1$. Hence
\begin{align*}
I+ II &\le c \intl_{t}^{0} (-s)^{ -\mu } \intl_{\R^n }   | v(s)|^2  | x|^{ \beta -1} dxds
\\
 &\le c \intl_{t}^{0} (-s)^{\beta(1-\mu ) -1} ds \le c (-t)^{\beta (1-\mu )}.  
\end{align*}
In order to estimate $ III + IV$ we make use of \eqref{3.1h} with $ \gamma = \beta -2$ and apply Cauchy-Schwarz's and Young's inequality  to obtain 
\begin{align*}
III+ IV &\le  c\intl_{t}^{0}   \bigg( \intl_{\R^n } | p|^2 | x|^{ \beta-2 }  dx\bigg)^{ 1/2} 
\bigg( \intl_{\R^n } | v(s)|^2 \eta _R(| x|)^2  | x|^{ \beta } dx\bigg)^{ 1/2}  ds
\\
&\le  c \intl_{t}^{0}  (-s)^{\beta(\frac12 -\frac{\mu}{2} ) -1}   \bigg( \intl_{\R^n } | v(s)|^2 \eta _R(| x|)^2 | x|^{ \beta }  dx\bigg)^{ 1/2}  ds
\\
&\le  c (-t)^{ \beta(1-\mu )} + \frac{1}{4} \esssup_{s\in (t, 0)} \intl_{\R^n } | v(s)|^2 \eta _R(| x|)^2  | x|^{ \beta }  dx.
\end{align*} 

 Inserting the estimates of $ I,II, III$ and $ IV$ into the right hand side of \eqref{3.1b},  and passing $R\to +\infty$, we get the desired estimate 
 \eqref{3.1a}.
 \hfill \Beweisende

\subsection{Fast decay using the local pressure for  exterior domains}

Let $ 0< r < +\infty$ be fixed.   By $ B(r)$ we denote the usual ball in $ \R^n $ with radius $ r>0$ 
with respect to the Euclidian norm having its center at  the origin.
For notational convenience by $ E^{ \ast}_r$ we denote the 
projection $ E^{ \ast}_{ B(r)^c}$ in $ D^{ -1,2}(B(r)^c)$ onto the closed subspace containing   functionals of the form $ \nabla \pi $,  which has been introduced in Section\,2.  Recalling the definition   $ E_r^{ \ast}$, we see that for every functional $ f\in D^{ -1,2}(B(r)^c)$ there exists a unique $ \pi \in L^2(B(r)^c)$ such that $ E_r^{ \ast}(f) = \nabla \pi $.

\begin{lem}\label{decay}
Let $ v\in L^2(-1, 0; L^2_\sigma (\R^n ))\cap  L^\infty_{ loc}([-1,0), W^{1,\, \infty}(\R^n ))$ be a  solution to the Euler equations  \eqref{1.1} satisfying \eqref{rate1}  for some $\mu \in [\frac{n}{n+2},1)$ and \eqref{conc2} with $ \sigma _0 = E  \delta _0$. Then for all  
$ k\in \N\cup \{ 0\}$ and $ 0< r <+\infty$ it holds 
\begin{equation}
\| v(t)- E^{ \ast}_{ r } (v(t))\|^2_{ L^2( B(r)^c)}
\le C_0^k 4^{ k^2}(-t)^{ (1-\mu )k} r^{ -k}\quad \forall\,t\in (-1,0),
\label{3.2}
\end{equation}
where the constant $ C_0>0$ depends only on $C_1$  of \eqref{1.3d} and $\mu$. 
\end{lem}

In the proof of Lemma\,\ref{decay} we make use of the following pressure estimate 

\begin{lem}
\label{lem3.7}
Let $ \Omega \subset \R^{n}$ be an exterior domain. Let $ \pi  \in L^q(\Omega )$  and $ f\in L^q(\Omega; \R^{n^2} )$, 
$ \frac{n}{n-1} < q < n$,  
such that $ \Delta \pi = \nabla \cdot \nabla \cdot f $ in $ \Omega $ in the sense of distributions.  Furthermore, let $ \zeta \in 
C^{\infty}(\Omega )$ such that $ \nabla \zeta \in C^{\infty}_c(\Omega )$.  Then   
\begin{align}
\| \pi \zeta^{ n} \|_{ L^q(K)} 
&\le c \| f\zeta^n \|_{ L^q} + 
c \Big(\max | \nabla \zeta | +  \mes(K)^{ 1/n}
\max (| \nabla^2 \zeta |+ | \nabla \zeta |^2) \Big)\times 
\cr
&\qquad \qquad \times \Big(\| f \zeta^{ n-2} \|_{ L^{ \frac{nq}{n+q}}(K)} + \| \pi \zeta ^{ n-2}\|_
{ L^{ \frac{nq}{n+q}}(K)}\Big),
\label{A.4a}
\end{align}
with a constant $ c>0$ depending only on $ n$ and $ q$, where $ K= \supp(\nabla \zeta )$. 
\end{lem}

{\bf Proof}: In our discussion below   we use the convention that repeated indices imply summation from $ 1$ to $ n$. By straightforward calculation we find that 
\begin{align*}
\Delta (\pi \zeta^{ n} ) &= \nabla \cdot \nabla \cdot (f \zeta^n)   + f: \nabla ^2\zeta^n -   
\partial _i ( f_{ ij} \partial _j\zeta^n) - \partial _j ( f_{ ij} \partial _i\zeta^n)  
\\
&\qquad \qquad \qquad + 2 \nabla \cdot  (\pi \cdot \nabla \zeta^n)  - \pi \Delta \zeta^n.  
\\
&= G_1+ G_2+ G_3+ G_4 + G_5 + G_6. 
\end{align*}
We may decompose $ \pi  \zeta^n $ into the sum $ \pi _1+\pi _2+ \pi _3+\pi _4+\pi _5+ \pi _6$,  where 
$$
\pi _i = N \ast G_i,\quad  i=1, \ldots, 6,
$$
and $N=N(x)$ is the fundamental solution of the Laplace equation in $\Bbb R^n$, given by
$$
 N(x) =\left\{\aligned & \frac{1}{c_n | x|^{ n-2}}\quad \text{for}\quad  n \ge 3,\\
 &-  \frac{1}{2\pi }\log(| x|) \quad \text{ for}\quad   n=2. 
 \endaligned
 \right.
$$

Using the Calder\'on-Zygmund estimate, we get 
\begin{align*}
\| \pi_1 \|_{ L^q} &\le c \| f \zeta^n \|_{ L^q},
\\
\|  \pi _3\|_{ L^q} + \| \pi _4\|_{ L^q}  + \| \pi _5\|_{ L^q}  &\le c\max | \nabla \zeta |
\Big(\| f \zeta ^{ n-1} \|_{ L^{ \frac{nq}{n+q}}(K)} + \| \pi \zeta ^{ n-1} \|_{ L^{ \frac{nq}{n+q}}(K)}\Big),
\\
\|  \pi _2\|_{ L^{ \frac{nq}{n-q}}} + \| \pi _6\|_{ L^{ \frac{nq}{n-q}}}  
 &\le c
\max (| \nabla^2 \zeta |+ | \nabla \zeta |^2)\times 
\\
&\qquad \times \Big(\| f \zeta ^{ n-2}\|_{ L^{ \frac{nq}{n+q}}(K)} + \| \pi \zeta ^{ n-2}\|_{ L^{ \frac{nq}{n+q}}(K)}\Big).
\end{align*}
Applying Jensen's inequality, we find 
\begin{align*}
&\|  \pi _2\|_{ L^{q}(K)} +\|  \pi _6\|_{ L^{q}(K)}
\\
&\qquad \le c \mes(K)^{ 1/n}
\max (| \nabla^2 \zeta | + | \nabla\zeta |^2 )
\Big(\| f \zeta ^{ n-2}\|_{ L^{ \frac{nq}{n+q}}(K)} + \| \pi \zeta ^{ n-2} \|_{ L^{ \frac{nq}{n+q}}(K)}\Big).
\end{align*}
Using triangle inequality together with the estimates of $ \pi _i$, $ i=1, \ldots, 6$, we obtain \eqref{A.4a}.  \hfill \Beweisende 

\begin{rem}
\label{rem3.8}
We may apply  Lemma\,\ref{lem3.7} for the case $q=2$ and $ \Omega = B(r)^{ c}$. If $ \zeta \in C^{\infty}(\R^n )$ is a cut off function 
such that $ \zeta \equiv 0 $ in $ B(2r)$, $ \zeta \equiv 1$ on $ B(4r)^c$ and $ | \nabla \zeta |^2 + | \nabla ^2 \zeta | \le c r^{ -2}$. 
Then the estimate \eqref{A.4a}  becomes 
 \begin{equation}
\| \pi \zeta^n \|_{L^2(K)} \le c \| f\zeta^n \|_{ L^2} + 
c r^{ -1}\Big(\| f \zeta ^{ n-2}\|_{ L^{\frac{2n}{n+2}}(K)} + \| \pi \zeta ^{ n-2}\|_{ L^{ \frac{2n}{n+2}}(K)}\Big). 
\label{A.4b}
\end{equation} 
Using H\"older's inequality and Young's inequality, we deduce from \eqref{A.4b}
\begin{equation}
\| \pi \zeta^n \|_{L^2(K)} \le c \| f\zeta^n \|_{ L^2} + 
c r^{ -\frac{n-1}{2}}\Big(\| f \|_{ L^{\frac{2n}{2n-1}}(K)} + \| \pi \|_{ L^{ \frac{2n}{2n-1}}(K)}\Big). 
\label{A.4c}
\end{equation}

\end{rem}

{\bf Proof of Lemma\,\ref{decay}}:  We prove  \eqref{3.2} by induction. Thanks to \eqref{2.8} having   
$$ \| E^{ \ast}_r (v(s))\|_{ L^2(B(r)^c)} \le  c\| v(s)\|_{ L^2} = c\| v(-1)\|_{ L^2},$$  the assertion is true for $ k=0$. 

\hspace{0.5cm}
We now assume \eqref{3.2}  is true for $ k\in \N\cup \{ 0\}$. 
Let $ 0<  r < +\infty$ be arbitrarily chosen, but fixed. In case $ 0 < r \le  4(-t)^{ 1-\mu }$ the  assertion is trivially fulfilled.  This can be readily seen 
by 
\[
\| v(t)- E^{ \ast}_{ r } (v(t))\|_{ L^2( B(r)^c)}^2 \le c \| v(t)\|^2_{ L^2} \le 
c\| v(-1)\|^2_{ L^2}\Big((-t)^{ 1-\mu } r^{-1}\Big)^k.
\]
Thus we only need to prove \eqref{3.2} for the opposite case  
\begin{equation}
\label{r4}
r > 4(-t)^{ 1-\mu }.
 \end{equation}
For notational simplicity we set 
\[
U:= B(r/4)^c, \quad  U_1:= B(r)^c, \quad  U_2:= B(r/2)^c.  
\]
Let $ \zeta \in C^{\infty}(\R^n )$ denote a cut off function such that 
$ 0 \le  \zeta \le  1$ in $ \R^n $, $ \zeta \equiv  0$ in $ B(r/2)$ and 
$ \zeta  \equiv 1$ on $U_1$.   As in Section 2 we define 
\begin{align*}
 \nabla p_{ h} &= - E^{ \ast}_{r /4}(v),\quad  \nabla  p_0 = - E^{ \ast}_{ r/4}( \nabla \cdot (v \otimes  v) ), 
\\
&\quad  \widetilde{v}  = v + \nabla p_{ h} = v - E^{ \ast}_{r/4}(v).
\end{align*}
Note that according to  \eqref{3.1a} it holds $ v(s)\in L^{ \frac{2n}{n+2}}(U) \subset D^{ -1,2}(U )$, and thus 
\[
\| p_h(s)\|_{ L^2(U)} \le c \| v(s)\|_{ L^{ \frac{2n}{n+2}}(U)}\quad \forall\,s\in (-1,0). 
\]  
Consulting \cite[Theorem\,A.4]{wol2} (with $ X= D_{0}^{ 1,2}(U ),  E^{ \ast}= E^{ \ast}_r$),  we see that the restriction of $ \nabla p$  to $ U$ equals 
to $ \partial _t\nabla p_h + \nabla p_0$ in the sense of the distribution,  i.e. the following identity holds true for all $ \varphi \in C^{\infty}_{c}(U\times (-1, 0))$, 
\begin{align*}
-  \intl_{-1}^{0}  \intl_{U} v \cdot \partial _t \varphi  + v \otimes v : \nabla \varphi dx dt
 = \intl_{-1}^{0}  \intl_{U}  \nabla p_h \cdot \partial _t \varphi    + p_0 \nabla \cdot \varphi  dx dt.   
\end{align*}
This shows that $ \widetilde{v} $ is a  solution to 
\begin{equation}
\partial _t \widetilde{v} +  (v \cdot \nabla) v = - \nabla p_0 \quad  \text{ in}\quad  U \times (-1,0).
\label{3.2b}
\end{equation}

 We compute
\begin{align*}
(v\cdot  \nabla) v &= (v\cdot \nabla)  \widetilde{v} - (v\cdot \nabla)  \nabla p_h =(v\cdot \nabla)  \widetilde{v} - (  \widetilde{v} -  \nabla p_h)\cdot \nabla ^2 p_h \\
&= 
(v\cdot \nabla) \widetilde{v}  - \widetilde{v}\cdot  \nabla^2 p_h
+  \frac{1}{2} \nabla | \nabla p_h|^2.
\end{align*}
Hence,  \eqref{3.2b}  implies that $ \widetilde{v} $ is a solution to the following transformed Euler equations 
\begin{equation}
\partial _t \widetilde{v} +(v\cdot \nabla) \widetilde{v} = \widetilde{v}\cdot  \nabla^2 p_h
-\nabla  p_1\quad  \text{ in}\quad  U \times (-1,0), 
\label{3.3}
\end{equation}
where we set
\[
p_1 = \frac{| \nabla p_h|^2}{2} + p_0. 
\]
Observing that $E^{ \ast}_{r /4} (\partial _t\widetilde{v} ) =0$ in the sense of distribution, we get 
\[
\nabla  p_1 = E^{ \ast}_{ r/4} \Big( -\nabla \cdot (\widetilde{v} \otimes  \nabla p_h)\Big) + E^{ \ast}_{ r /4} \Big( \nabla \cdot ( v \otimes   \widetilde{v}) \Big)
=:\nabla  p_{ 11}+ \nabla  p_{ 12}.
\]
Let $ t \le  s < 0$ be fixed.  Since $ \Delta p_{ 11}(s) = -\nabla \cdot \nabla \cdot (\widetilde{v}(s) \otimes  \nabla p_h(s))$, using Lemma\,3.7 and Remark\,\ref{rem3.8}, we find 
 \begin{align*}
&\|p_{ 11}(s)\zeta^n \|_{ L^2(K)} 
\\
& \quad \le c \Big\{ \| (\widetilde{v}(s) \otimes  \nabla p_h(s)) \zeta ^n\|_{ L^2(U)}+ 
r^{ - \frac{n-1}{2}}\| \widetilde{v} \otimes  \nabla p_h(s) \|_{ L^{ \frac{2n}{2n-1}}(K)} 
\\
& \qquad \qquad \qquad \qquad \qquad + 
r^{ - \frac{n-1}{2}}\| p_{ 11}(s)\|_{ L^{ \frac{2n}{2n-1}}(K)}\Big\}
\\
&\quad  \le c \Big\{ \| (\widetilde{v}(s) \otimes  \nabla p_h(s)) \zeta ^n\|_{ L^2(U)}+ 
r^{ - \frac{n-1}{2}}\| \widetilde{v}(s) \otimes  \nabla p_h(s) \|_{ L^{ \frac{2n}{2n-1}}(U)}\Big\}, 
\end{align*}
where $ K = \supp (\nabla \zeta )$.  Applying H\"older's inequality,  we infer 
\begin{align*}
\| \widetilde{v} \otimes  \nabla p_h(s) \|_{ L^{ \frac{2n}{2n-1}}(U)} &\le c\| \widetilde{v} (s)\|_{ L^2(U)} \| \nabla p_h(s)\|_{ L^{ \frac{2n}{n-1}}(U)}
\\
&\le c\| \widetilde{v} (s)\|_{ L^2(U)} \| v(s)\|_{ L^{ \frac{2n}{n-1}}(U)}
\\
&\le c  \| \widetilde{v} (s)\|_{ L^2(U)} \| v(s)\|_{ L^2}^{ \frac{n-1}{n}} \| v(s)\|_{ L^\infty}^{ \frac1n}.
\end{align*}
Furthermore,  since $\mu \ge \dfrac{n}{n+2}$,  we have $ \mu \dfrac{n-1}{n} \ge  (1-\mu ) \dfrac{n-1}{2}$, and therefore from \eqref{r4} we obtain
\begin{align}
\| v(s)\|_{ L^\infty}^{ \frac1n} &\le C_1^{\frac1n}(-s)^{ - \frac{\a}{n}}  = 
C_1^{\frac1n}(-t)^{ \mu \frac{n-1}{n}} (-s)^{ -\mu } \cr 
&\label{elem}\le C_1^{\frac1n} (-t)^{(1- \mu) \frac{n-1}{2} }(-s)^{ -\mu } 
\le cr^{ \frac{n-1}{2}} (-s)^{ -\mu }.
\end{align}
Hence, we estimate 
\[
 r^{ - \frac{n-1}{2}}\| \widetilde{v} \otimes  \nabla p_h(s) \|_{ L^{ \frac{2n}{2n-1}}(U)} \le  c  C_1 \| \widetilde{v} (s)\|_{ L^2(U)}   (-s)^{ -\mu }. 
\]

Similarly,  
\begin{align*}
\| \widetilde{v}(s) \otimes  \nabla p_h(s)\zeta^n \|_{ L^2(U)} &\le c r^{ - \frac{n-1}{2}} \| \widetilde{v}(s) \|_{ L^2(U)} 
\| \nabla p_h(s)\|_{ L^{ \frac{2n}{n-1}}(U)}
\\
&\le c  r^{ - \frac{n-1}{2}}\| \widetilde{v} (s)\|_{ L^2(U)} \| v(s)\|_{ L^2}^{ \frac{n-1}{n}} \| v(s)\|_{ L^\infty}^{ \frac1n}
\\
&\le c  C_1 \| \widetilde{v} (s)\|_{ L^2(U)}   (-s)^{ -\mu }. 
\end{align*}
This shows that 
\begin{align*}
\|p_{ 11}(s)\zeta^n \|_{ L^2(K)}  \le cC_1 \| \widetilde{v} (s)\|_{ L^2(U)}   (-s)^{ -\mu }. 
\end{align*}
Similarly  we get 
\[
\|p_{ 12}(s)\zeta^n \|_{ L^2(K)} \le   c C_1\| \widetilde{v} (s)\|_{ L^2(U)}   (-s)^{ -\mu }. 
\]
  We now assume \eqref{3.2} is true for $\tilde{v}(t)= v(t)- E^{ \ast}_{ r/4 } (v(t))$ with  $k$, then inserting this into the above estimates for $p_{ 11}, p_{12}$,  we find 
\begin{equation}
\| p_1(s) \zeta^n \|_{ L^2(K)} 
\le c  C_1^{\frac12}  C_0^{ \frac{k}{2}}  4^{\frac{k^2}{2}+ k}  (-s)^{ \frac12(1-\mu )- \mu }  r^{-\frac{k}{2}}. 
\label{3.4}
\end{equation}

\hspace{0.5cm}
We multiply  \eqref{3.3}  by $ -\widetilde{v} \zeta^{ 2n} $,  integrate the result over $ U \times (t,0), -1<t<0$, and apply integration by parts. This yields
\begin{align}
&\frac{1}{2} \| \widetilde{v}(t) \zeta ^n\|_{ L^2(U)}^2 
\cr
&\qquad = -  \intl_{t}^{0} \intl_{U}  v(s)\cdot \nabla \zeta  | \widetilde{v}(s) |^2   \zeta^{ 2n-1} dxds - 
 \intl_{t}^{0} \intl_{U}  \widetilde{v}(s) \otimes \widetilde{v}(s) :\nabla ^2 p_h(s)    \zeta^{ 2n} dx ds
\cr
&\qquad \qquad +   2\intl_{t}^{0} \intl_{U} p_1(s) \widetilde{v}(s) \cdot \zeta^{ 2n-1} \nabla \zeta dx ds=I+II+III.
\label{3.5}
\end{align}
Applying Cauchy-Schwarz's inequality, and again using the assumption  of \eqref{3.2} for $k$ and $ \frac{r}{4}$ in place of $ r$ together with $ \| v(s)\|_{ L^\infty} \le c(-s)^{ (1-\mu )}$, we find  
\begin{align*}
I &\le cr^{ -1}\intl_{t}^{0}\| \widetilde{v} (s)\|^2_{ L^2(U)} \| v(s)\|_{ L^\infty}  ds\le c  r^{ -1}
4^{ k^2}C_0^k \Big(\frac{r}{4}\Big)^{ -k} \intl_{t}^{0} (-s)^{(1-\mu )k - \mu } ds 
\\
& \le  c C_1 C_0^k 4^{ (k+1)^2} (-t)^{ (1-\mu )(k+1)} r^{ -k-1}. 
\end{align*}
 Using Lemma\,\ref{lemA1a}, we estimate  
 \begin{align*}
 \| \nabla^2 p_h (s)\zeta^{ 2n}  \|_{ L^\infty(U)} &\le c r^{- \frac{n+1}{2}}\| \nabla p_h\|_{ L^{ \frac{2n}{n-1}}(U)} \le 
 c r^{- \frac{n+1}{2}} \| v(s)\|_{ L^{ \frac{2n}{n-1}}(U)} \\
 &\le c   r^{- \frac{n+1}{2}} \| v(s)\|_{ L^2}^{ \frac{n-1}{n}} \| v(s)\|_{ L^\infty}^{ \frac1n}
 \le c r^{ -1} (-s)^{ -\mu },  
 \end{align*}
where for the second inequality we have applied  \eqref{2.8} with $ q= \frac{2n}{n-1}$, while for the fourth inequality we  have used  \eqref{elem}. Thus, 
by similar reasoning as we have used for the estimation of $ I$,  we get 
\begin{align*}
II &\le  \intl_{t}^{0} \| \widetilde{v} (s)\|^ 2_{ L^2(U)} \| \nabla^2 p_h (s)\zeta^{ 2n}  \|_{ L^\infty(U)} ds  
\\
&\le 
c C_1 C_0^k 4^{k^2}    \Big(\frac{r}{4}\Big)^{ -k}  r^{ -1} \intl_{t}^{0}  (-s)^{  (1-\mu )k- \mu } ds
\\
& \le c C_1 C_0^k 4^{(k+1)^2}  (-t)^{  (1-\mu )(k+1)} r^{ -1-k}. 
\end{align*}
 Finally, applying Cauchy-Schwarz's inequality together  with \eqref{3.4}, and the assumption \eqref{3.2} for $ k$,  we estimate 
 \[
III \le  c C_1C_0^k 4^{ (k+1)^2}(-t)^{ (1-\mu )(k+1)} r^{ -k-1}. 
\] 
 Inserting the estimates of $ I, II$ and $ III$ into \eqref{3.5},  we are led to 
 \[
\| v(t) - E^{ \ast}_{ r /4}(v(t))\|^2_{ L^2(U_1)} \le 
c C_1 C_0^k 4^{ (k+1)^2}(-t)^{  (1-\mu )(k+1)} r^{ -k-1}.
\] 
Since 
\begin{align*}
 v(t) - E^{ \ast}_{ r}(v(t)) &= v(t) - E^{ \ast}_{ r/4 }(v(t)) + E^{ \ast}_{ r /4}(v(t)) - E^{ \ast}_{ r}(v(t)) 
\\
&= v(t) - E^{ \ast}_{ r/4}(v(t)) + E^{ \ast}_{r }(E^{ \ast}_{ r /4}(v(t))) - E^{ \ast}_{ r}(v(t)) 
\\
&= v(t) - E^{ \ast}_{ r/4}(v(t)) - E^{ \ast}_{ r}\Big(v(t) - E^{ \ast}_{ r /4}(v(t))\Big)
\\
&= \widetilde{v}(t)  - E^{ \ast}_{ r} (\widetilde{v}(t))
\end{align*}
in $ U_1$, we estimate 
\begin{align*}
\| v(t) - E^{ \ast}_{ r}(v(t))\|^2_{ L^2(    B(r)^c)} &=\|\widetilde{v}(t)  - E^{ \ast}_{ r} (\widetilde{v}(t))\|^2_{ L^2(    B(r)^c)} \le c \|\widetilde{v}(t) 
\|^2_{ L^2(    B(r)^c)}\\
&= 
c\| v(t) - E^{ \ast}_{ r/4}(v(t))\|^2_{ L^2(U_1)} 
\\
&\le c C_1 C^k_0 4^{ (k+1)^2}(-t)^{ (1-\mu )(k+1)} r^{ -k-1}.
\end{align*} 
 This shows that \eqref{3.2} holds for $ k+1$ with $ C_0 = c C_1$.   
  \hfill \Beweisende

\subsection{Proof of Theorem\,\ref{thm3.1}}
  
Let  us fix $\theta$ so that 
\begin{equation}
\label{theta}
0<\theta < \frac{1}{n+2}.
\end{equation}
 For given solution  $(v, p)$ to the Euler equations we define 
 \begin{align*}
w(x,t) &= v ( (-t)^{ \theta }x, t  ),
  \\
\pi (x,t) &=(-t)^{ -\theta } p ( (-t)^{ \theta }x, t  ),\quad  (x,t) \in \R^{n}\times (-1,0).
  \end{align*}
Then,  $(w, \pi )$ solves
\begin{align}
\frac{\partial w}{\partial t}  + \theta (-t)^{ - 1} x\cdot \nabla w
 + (-t)^{ - \theta }(w\cdot \nabla) w &= - \nabla \pi,
 \label{6.3b}
 \\
 \nabla \cdot w &=0.
 \label{6.3a}
\end{align}  
Using the transformation formula, we find  
\[
\| w(t)\|_{ L^2}^2 = (-t)^{ -n\theta } \| v(t)\|^2_{ L^2} =  (-t)^{ -n\theta } \| v(-1)\|^2_{ L^2}. 
\]

On the other hand, by Lemma \ref{spacedecay} we infer that for any $ 0< \beta < n+2$
\begin{align*}
\| w(t)\|_{ L^2(B(1)^c)}^2 &= (-t)^{ -n\theta  } \| v(t)\|^2_{ L^2(B((-t)^{ \theta })^c)}\\
& \le (-t)^{ -n\theta  }\int_{\{|x|>(-t)^\theta\}} |v(t)|^2 \frac{|x|^\beta }{(-t)^{\beta\theta}} dx
\\
&\le C (-t)^{ -n\theta  } (-t)^{ -\beta\theta } 
(-t)^{\frac{2}{n+2} \beta } = C (-t)^{ -(n+\beta)\theta + \frac{2}{n+2}\beta  }.  
\end{align*}
Choosing  $\beta =n$,  we get $  -(n+\beta)\theta +\frac{2}{n+2}\beta >0$ for $\theta$  satisfying \eqref{theta}.  Therefore 
\begin{equation}
\label{wconv}
\lim_{t\to 0} \| w(t)\|_{ L^2(B(1)^c)}^2 =0.
\end{equation}

\hspace{0.5cm}
By $ \PP_r$, $ 0< r< +\infty$ we denote the Helmholtz projection from $ L^2(B(r)^c)$ onto $ L^2_{ \sigma }(B(r)^c)$. We easily calculate 
\[
 \PP_1 w(x, t) = (\PP_{ (-t)^\theta} v)((-t)^{ \theta }x, t). 
\]
To see this we only need to check that $ w(x,t) - (\PP_{ (-t)^\theta} v)((-t)^{ \theta }x, t)$ is a gradient  field. Indeed, 
\begin{align*}
 w(x,t) - (\PP_{ (-t)^\theta} v)((-t)^{ \theta }x, t) &= v((-t)^{ \theta }x, t) - (\PP_{ (-t)^\theta} v)((-t)^{ \theta }x, t)
\\
&=
 \nabla q((-t)^{ \theta }x, t).
\end{align*}
Appealing to Lemma\,\ref{decay} for $ \mu = \frac{n}{n+2}$, we see that for  every $ r>0$ and $ k\in \N$ it holds 
\[
\| v(t) - E^{ \ast}_{ r} (v(t))\|_{ L^2(B(r)^c)}^2 \le C(k)  (-t)^{ \frac{2}{n+2}k} r^{ -k},
\]
where $ C(k)$ depends on $ k$ and  $C_1$ only. Noting that 
\[
 \PP_{ r} v(t) = \PP_{ r}  (v(t) - E^{ \ast}_{ r} (v(t))), 
\]
 from the above estimate we deduce 
\[
\|\PP_r v(t)\|_{ L^2(B(r)^c)}^2 \le   \| (v(t) - E^{ \ast}_{ r} (v(t)) \|_{ L^2(B(r)^c)}^2  \le C(k)   (-t)^{\frac{2}{n+2}k} r^{ -k}.
\]
This yields 
\begin{align*}
\| \PP_1 w(t)\|_{ L^2(B(1)^c)}^2 &\le (-t)^{ -n\theta  }  \| \PP_{ (-t)^{ \theta }} (v(t))\|_{ L^2(B((-t)^\theta )^c)}^2
\\
&\le  C(k) (-t)^{ \frac{2}{n+2}k} (-t)^{ -n\theta  } (-t)^{ -k\theta } = C(k) (-t)^{k  (\frac{2}{n+2}- \theta ) -n\theta }.
\end{align*}
Since $\theta $ satisfies \eqref{theta},  this shows the decay rate of  $ \| \PP_{ 1} (w(t))\|_{ L^2(B(1)^c)}^2$ as $ t \rightarrow 0$  is of any order $ O((-t)^k)$.  

\hspace{0.5cm}
Now we set $ w_0(t) := \PP_1(w(t))$ on $ B(1)^c$ and $ \nabla q_h(t) = w(t) - w_0(t)$.  Since  $ \nabla \cdot w=0$,  we see that $ \nabla q_h(t)$ 
is harmonic,  and therefore it also solves the system \eqref{6.3b}-\eqref{6.3a} with 
\[
\widetilde{\pi} = \partial _t q_h + \theta (-t)^{ -1} \Big(x\cdot \nabla q_h - q_h+ \frac{1}{2} | \nabla q_h|^2\Big)
\]
in place of $ \pi $.  Taking the difference of the two equations for $w$ and $\nabla q_h$ respectively,    we get 
\begin{align}
\frac{\partial w_0}{\partial t}  + \theta (-t)^{ - 1} x\cdot \nabla w_0
 + (-t)^{ - \theta }(w\cdot \nabla)  w-  (-t)^{ - \theta }\nabla q_h\cdot \nabla^2 q_h &= - \nabla (\pi -\tilde{ \pi}),
 \label{6.6}\\
 \nabla \cdot w_0 &=0,
 \label{6.5}
\end{align}  
the both of which  are in $ B(1)^c \times (-1, 0)$.  Note that 
\begin{align*}
(w\cdot \nabla)  w-  \nabla q_h\cdot \nabla^2 q_h &=
(w_0\cdot \nabla)   w + (\nabla q_h\cdot \nabla) w -  \nabla q_h\cdot \nabla^2 q_h 
\\
&= (w_0\cdot \nabla)  w + (\nabla q_h\cdot \nabla) w_0 .
\end{align*}
Therefore,  \eqref{6.6} turns into  
\begin{equation}
\frac{\partial w_0}{\partial t}  + \theta (-t)^{ - 1} x\cdot \nabla w_0
 + (-t)^{ - \theta }(w_0\cdot \nabla)  w+  (-t)^{ - \theta }(\nabla q_h\cdot \nabla) w_0 = - \nabla (\pi -\tilde{ \pi}).  
\label{6.7}
\end{equation}
We now multiply \eqref{6.7} by $- w_0(s)$, integrate it over $ B(1)^c \times (t,0)$,  and then apply the integration by parts. Taking into account \eqref{wconv}, we have
the identity 
\begin{align}
&\frac{1}{2} \| w_0(t)\|_{ L^2(B(1)^c)}^2 + \frac{n\theta  }{2}  \intl_{t}^{0} \intl_{B(1)^c} (-s)^{ -1} | w_0(s)|^2 dx ds
+  \frac{\theta }{2} \intl_{t}^{0} \intl_{\partial B(1)} (-s)^{ -1} | w_0|^2 dS ds\cr
&=  \intl_{t}^{0} \intl_{B(1)^c}  (-s)^{ -\theta }w_0(s) \otimes w_0(s): \nabla w(s) dxds - \frac{1}{2}\intl_{t}^{0} \intl_{\partial B(1)}  (-s)^{ -\theta }  w(s)\cdot x | w_0(s)|^2  dSds,
\label{www}\end{align}
where we used the fact 
$$
\intl_{B(1)^c} \nabla q_h(s) \cdot \nabla | w_0(s)|^2 dx = \intl_{B(1)^c} w(s) \cdot \nabla | w_0(s)|^2 dx
=   -\intl_{\partial B(1)} x\cdot w (s) | w_0(s)|^2dS
$$
for the second integral of the  right-hand side.
Since $-\theta-\frac{n}{n+2} >-1$ due to \eqref{theta}, we may choose $ -1< t_0 < 0$ so  that 
\begin{align*}
\max_{ x\in \partial B(1)} (-s)^{ -\theta }| w(x, s)| &\le  (-s)^{ -\theta } \|v(s)\|_{L^\infty} \le c (-s)^{-\theta-\frac{n}{n+2} } \\
& \le  \frac{\theta }{2} (-s)^{ -1}\quad  \forall\,t_0 \le  s < 0,
\end{align*}
which implies that the second term of the  right-hand side of \eqref{www} can be absorbed into the third term of the 
 left-hand side of \eqref{www}.
Then, since $ (-s)^{ -\theta } | \nabla w(s)| \le \|\nabla v(s)\|_{L^\infty} \le  \frac{a}{2} (-s)^{ -1} $ for all $s\in (-1, 0)$, where we set
$$a=2\sup_{-1<s<0} (-s) \|\nabla v(s)\|_{L^\infty},$$

 we obtain 
\begin{equation}
\| w_0(t)\|_{ L^2(B(1)^c)}^2 \le a \intl_{t}^{0} \intl_{B(1)^c}  (-s)^{ -1}| w_0(s)|^2  dxds\qquad \forall t\in (t_0, 0).  
\label{6.8}
\end{equation}
Let us define 
\[
X(t) :=  \intl_{t}^{0} \intl_{B(1)^c}  (-s)^{ -1}| w_0(s)|^2  dxds. 
\]
Then, from \eqref{6.8} it follows that 
\begin{equation}
- (-t) X'(t)  \le a X(t)\qquad \forall t\in (t_0, 0),\label{6.9}
\end{equation}
which  is equivalent to $ X'(t) \ge - a(-t)^{ -1} X(t)$. If we assume that $ X(t) >0$ for all $ t\in (t_0, 0)$, we get 
\[
(\log X)' \ge  (\log (-t)^a)'\quad  \Leftrightarrow \quad  \Big(\log \frac{X(t)}{(-t)^a}\Big)' \ge 0. 
\]
Accordingly, $\log \frac{X(t)}{(-t)^a} $ is nondecreasing, and by the monotonicity of $ \log$ the function 
 $ \frac{X(t)}{(-t)^a}$ is also nondecreasing. However, by the fast decay of $ \| w_0(t)\|_{ L^2(B(1)^c)}$ as $ t \rightarrow 0$ 
 $ X(t)$ is decaying faster to zero than $ (-t)^a$. Therefore 
 \[
\lim_{t \to 0} \frac{X(t)}{(-t)^a}=0,
\]
which  is a contradiction to $ \frac{X(t)}{(-t)^a} \ge \frac{X(t_0)}{(-t_0)^a} >0$  for all $ t\in (t_0,0)$. Consequently, 
$ X \equiv 0$. This shows that $ \nabla \times w(t) =0 $ on $ B(1)^c$ for all $ t_0< t< 0$.  This implies that the  vorticity $\omega (t)=\nabla \times v(t)$ also vanishes on $ B((-t)^\theta )^c$ for all $ t_0 < t< 0$, namely
$$
\{ x\in \Bbb \R^n \, |\, |\omega (x,t) |>0\} \subset B((-t)^\theta )\quad \forall t\in (t_0, 0).
$$
Since the measure of the set $ \{x\in \R^n  \,|\, | \omega (x, t)| >0\}$ is conserved for $t\in (-1, 0)$  by virtue of the vorticity transport formula (see e.g.\cite[Proposition 1.8]{maj}), we have
\begin{align*}
 \mes \{x\in \R^n  \,|\, | \omega (x, t_0)| >0\} &= \mes 
  \{x\in \R^n  \,|\, | \omega (x, t)| >0\} \le c(-t)^{ n\theta  } 
\\
&\rightarrow 0 \quad  \text{ as}\quad  t \rightarrow 0. 
\end{align*}
Whence, $ \omega(t_0) \equiv 0$,  which implies that  $ v(t_0)$ is harmonic. Recalling that 
$ v(t_0)\in  L^2(\R^n )$, we conclude that $ v(t_0) \equiv 0$, and hence  $ v \equiv 0$. 
 \hfill \Beweisende  

\section{Proof of Theorem\,\ref{thm1.4}}
\label{sec:-?}
\setcounter{secnum}{\value{section} \setcounter{equation}{0}
\renewcommand{\theequation}{\mbox{\arabic{secnum}.\arabic{equation}}}}

\subsection{Local criterion for the energy non-concentration}
\label{sec:-2}
\setcounter{secnum}{\value{section} \setcounter{equation}{0}
\renewcommand{\theequation}{\mbox{\arabic{secnum}.\arabic{equation}}}}

In  this first subsection we remove one point energy concentration for local  
weak solution to the Euler equations  satisfying the  local energy inequality under a weaker condition  than the one in  Shvydkoy \cite{shv}.    
In our discussion below we make use of the following notation. 
We define the following space time cylinder  
\[
Q(r) = B(r) \times I(r),\quad  \text{ where} \quad  I(r) = (-r^{ \frac{n+2}{2}}, 0).
\]
Let $ 0< R< +\infty$ be fixed.  We consider the Euler equations   
\begin{equation}
\partial _t v+ (v\cdot \nabla) v = -\nabla p,\quad  \nabla \cdot v =0\quad  \text{ in}\quad  Q(R). 
\label{9.1}
\end{equation}
The main result of this subsection is the following 

\begin{thm}
\label{thm9.1}
Let $v \in L^\infty(I (R); L^2(B(R)))\cap L^3(Q(R))$ 
be  a local suitable weak solution to \eqref{9.1} according to Definition\,\ref{def2.1} 
such that the local energy inequality \eqref{2.10} is fulfilled.
Furthermore, we assume that
\begin{equation}
\sup_{0< r \le  R} r^{ -1} \| v\|_{ L^3(Q(r))}^3 < +\infty,\quad  \liminf_{r \to 0^+} r^{ -1} \| v\|_{ L^3(Q(r))}^3 =0. 
\label{9.4}
\end{equation}
Then, there is no energy concentration at the point  $ x=0$ as $t\to0^-$.  More precisely, if $ \sigma_0 \in \mathcal{M}_{ | v|^2}(0)$ then 
\begin{equation}
\sigma_0 (\{ 0\})=0. 
\label{conc9}
\end{equation}

\end{thm}
 \begin{rem} In \cite{shv} Shvydkoy showed that if $v\in L^{q} (-1, 0; L^\infty (\Omega))\cap L^\infty(-1, 0; L^2 (\Omega))$, $q=  \frac{2}{n}+1$,
 is a suitable weak solution, then  the measure in $\mathcal{M}_{ | v|^2}(0)$ has no atoms in $\Omega$.   This actually follows from the above theorem immediately.
 Indeed, let $Q(r)\subset \Omega \times (-1, 0)$, then
 \bqn
  r^{-1} \|v\|^3 _{L^3 (Q(r))} &=&r^{-1} \int_{-r^{\frac{n+2}{2}}} ^0 \int_{B(r)} |v|^3 dx dt \\
  &\leq& \|v\|_{L^\infty (-1,0; L^2 (\Omega))} ^2 r^{-1} 
 \int_{-r^{\frac{n+2}{2}}} ^0  \|v\|_{L^\infty(B(r))} dt \\
 & \leq& \|v\|_{L^\infty (-1,0; L^2 (\Omega))} ^2\left( \int_{-r^{\frac{n+2}{2}}} ^0  \|v\|_{L^\infty(B(r))}^{\frac{2}{n}+1} dt\right)^{\frac{n}{n+2}} \to 0
 \eqn
 as $r\to 0$.
\end{rem}
{\bf Proof of Theorem\,\ref{thm9.1}}:  Let $ \nabla p_h$ denote the local presssure $ - E^{ \ast}_{ B(R)}(v)$,  which has been defined in Definition\,\ref{def2.1}. 
By virtue of Lemma\,\ref{lem2.3} there exists a unique measure valued trace $ \widetilde{\sigma }\in L^\infty(I(R); \mathcal{M}^+ (B(R))) $ 
of the function $ | \widetilde{v} (\cdot )|^2$, where $ \widetilde{v} = v+ \nabla p_h$ (cf. also Remark\,\ref{rem2.4}).  
Thanks to Lemma\,\ref{lem2.4} we only need to show that $\widetilde{\sigma }_0:=  \widetilde{\sigma } (0)$ has no atoms. In fact, it suffices 
to prove that $ \widetilde{\sigma}  _0(\{ 0\})=0$.

\hspace{0.5cm}
  Following the arguments  of Section\,2,  we define the local pressure 
\[
\nabla p_h = -E^{ \ast}_{ B(R)}(v),\quad \nabla p_0 = -E^{ \ast}_{ B(R)}((v\cdot \nabla) v). 
\]

Recalling Definition\,\ref{def2.1}, the function
$ \widetilde{v} = v+\nabla p_h$ solves the equation 
\begin{equation}
\partial _t \widetilde{v} + (v\cdot \nabla) v = - \nabla p_0\quad  \text{ in}\quad  Q(R).
\label{9.5a}
\end{equation}
Since $ v$ is a local suitable weak solution to the  Euler equations  (cf. Definition\,\ref{def2.1}) by means  of  Lemma\,\ref{lem2.3}  the generalized local energy inequality is satisfied  
for a.e.  $t\in I(R)$, and for all nonnegative $ \phi \in C^\infty_c(B(R))$   
\begin{align}
&\intl_{B(R) } | \widetilde{v} (t)|^2 \phi   dx 
\cr
&\qquad \le   \intl_{B(R)}  \phi d \widetilde{\sigma_0}    	
+  \intl_{t}^{0} \intl_{B(R)}   |\widetilde{v}  |^2  v\cdot \nabla \phi  dx d\tau + 
 \intl_{t}^{0} \intl_{B(R)}  2 p_0  \widetilde{v} \cdot \nabla \phi  dx d\tau 
\cr
&\qquad \qquad + 
 \intl_{t}^{0} \intl_{B(R)}  v \otimes v: \nabla ^2 p_h  \phi  dx d\tau.
\label{9.5d}
\end{align}

Let $ 0 < r< \frac{R}{2}$.  Let $ \zeta \in C^{\infty}_{c}(B(R))$ denote a cut off function 
such that $ 0 \le \zeta  \le 1$ in $ B(R)$, $ \zeta \equiv 1$  on $ B(R/2)$. Furthermore, let $ \eta\in  C^{\infty}_{c}(B(r/2))$ denote  
a cut off function such that $ 0 \le \eta  \le 1$ in $ B(r/2)$, $ \eta \equiv 1$  on $ B(r/4)$, and $ | \nabla^k \eta | \le c_k r^{ -k}$ for all 
$k\in \N $. Let 
$ - r^{ \frac{n+2}{2} }< t < 0$. In \eqref{9.5d} we 
put $ \phi= \zeta (1- \eta ) $.  This yields 
\begin{align}
&\intl_{ B(R)} | \widetilde{v} (t)|^2 \zeta (1-\eta ) dx 
\\
&\le  \intl_{B(R)} \zeta (1- \eta ) d\widetilde{\sigma} _0 
+    \intl_{t}^{0} \intl_{B(R)}   | \widetilde{v}  |^2  v\cdot \nabla \zeta   dx d\tau
- \intl_{t}^{0} \intl_{B(r)}   | \widetilde{v}  |^2  v\cdot \nabla \eta   dx d\tau
\cr
& \qquad+    \intl_{t}^{0} \intl_{B(R)}  2 p_0 \widetilde{v} \cdot \nabla \zeta   dx d\tau
- \intl_{t}^{0} \intl_{B(r)}   2 p_0 \widetilde{v} \cdot \nabla \eta   dx d\tau
\cr
&\qquad + \intl_{t}^{0} \intl_{B(R)}  v \otimes v: \nabla ^2 p_h  \zeta (1-\eta )  dx d\tau
 \cr
&\quad =  \intl_{B(R)} \zeta (1- \eta ) d\widetilde{\sigma} _0   +I+ II+ III+ IV+ V.
\label{9.7}
\end{align}

 In our discussion below we frequently make use of the following inequalities for almost every  $ \tau \in (-1,0)$
 \begin{align}
 \| p_0(\tau )\|_{ L^{ 3/2}(B(R))} &\le  
c\| v(\tau )\|_{ L^3(B(R))}^2,
\label{9.7a}   
\\
\| \nabla p_h(\tau )\|_{ L^3(B(R))} &\le c\| v(\tau )\|_{ L^3(B(r))} 
\label{9.7b}
\end{align}

By means of H\"older's inequality and  Young's inequality,  using \eqref{9.7b}, we easily get 
\[
I +II\le c R^{ -1} \| v\|^3_{ L^3(B(R)\times I(r))} + c r^{ -1} \| v\|^3_{ L^3(Q(r)} + c r^{ -1}  
\| \nabla p_h\|^3_{ L^3(Q(r)}
\]
Recalling that $ \nabla p_h(\tau ) $ is harmonic in $ B(R)$ and employing  \eqref{9.7b}, we get for the last term on the right-hand side of the above inequality 
\[
r^{ -1}  \| \nabla p_h\|^3_{ L^3(Q(r)} \le c R^{ -1}  \| \nabla p_h\|^3_{ L^3(B(R)\times I(r)}  
\le R^{ -1}  \| v\|^3_{ L^3(B(R)\times I(r)}. 
\]
Combining the last two inequalities, we arrive at 

\[
I + II\le c R^{ -1} \| v\|^3_{ L^3(B(R)\times I(r))} + c r^{ -1} \| v\|^3_{ L^3(Q(r))}. 
\]
Applying H\"older's inequality and using  \eqref{9.7a} 
for  almost every  $ \tau \in I(r)$, we get 
\[
III\le c R^{ -1} \| v\|^3_{ L^3(B(R)\times I(r))}. 
\] 

We proceed to estimate $ V$.  By virtue of  Sobolev's embedding  theorem  we see that  $ W^{m,\, 2}(B(R)) \hookrightarrow  L^3(B(R)) $ for $ m \ge \frac{n}{6}$. This  together with  Lemma\,\ref{lemA.1}  and \eqref{9.7b} gives
\begin{align*}
\| \nabla ^2 p_h(s) \zeta \|_{ L^3(B(R))} &\le  \sum_{k=0}^{m} R^{ -\frac{n}{6}+k}
 \| \nabla ^k(\nabla^{ 2} p_h(s) \zeta )\|_{ L^{2}(B(R))} 
\\
&\le c R^{ -\frac{n}{6} -1}\| \nabla  p_h (s)  \|_{ L^2(B(R))} \le  c R^{- \frac{n}{6} -1}\|v(s) \|_{ L^2(B(R)))}
\\
 &\le  cR^{ -1}\|v(s) \|_{ L^3(B(R)))}. 
\end{align*}
Using  H\"older's inequality together with the above estimate of $ \nabla  ^2 p_h$ we obtain
\[
 V \le c R^{ -1} \| v\|^3_{ L^3(B(R)\times I(r))}. 
\]

\hspace{0.5cm}
It only  remains the estimate the integral $ IV$,  which contains the pressure $ p_0$.  Observing the condition \eqref{9.4}, we find that 
\[
 \sup_{ 0< \rho  \le R}  \rho ^{ -1}\| v \otimes v\|^{ 3/2}_{ L^{ 3/2}(Q(\rho ))} < +\infty. 
\] 

Applying Lemma\,\ref{lemB.1} with $ f= v \otimes v$, $ p= \frac{3}{2}$, and $ \lambda =1$ (cf. also \cite[Lemma 2.8]{cw1}), it can be checked that 
\begin{equation}
 \sup_{ 0< \rho  \le R} \rho ^{ -1}\| p_0\|^{ 3/2}_{ L^{ 3/2}(Q(\rho ))} < +\infty. 
\label{9.8}
\end{equation}
Applying H\"older's inequality along with \eqref{9.8}, we infer 
\[
IV \le c\sup_{ 0< \rho  \le R} \rho ^{ -1}\| p_0\|^{ 3/2}_{ L^{ 3/2}(Q(\rho ))}  (r^{ -1} \| v\|^3_{ L^3(Q(r))})^{ 1/3}. 
\]
Inserting the estimates of $ I, \ldots, VII$ into the right-hand side of \eqref{9.7}, we arrive at 
 \begin{align}
& \intl_{ B(R)} | \widetilde{v} (t)|^2 \zeta dx 
\cr
&\quad  \le  \intl_{B(R)} \zeta (1- \eta ) d\widetilde{\sigma} _0  +\intl_{ B(r)} | \widetilde{v} (t)|^2 \eta dx 
\cr
&\qquad \qquad  + c R^{ -1} \| v\|^3_{ L^3(B(R)\times I(r))}
 + 
c \Big(\frac{r}{R}\Big)^{ n/2} \| v\|^3_{ L^\infty(-R^{ 5/2}, 0; L^2(B(R)))} 
\cr
&\qquad \qquad + c r^{ -1} \| v\|^3_{ L^3(Q(r)}+ 
c\sup_{ 0< \rho  \le R} \rho ^{ -1}\| p_0\|^{ 3/2}_{ L^{ 3/2}(Q(\rho ))}  (r^{ -1} \| v\|^3_{ L^3(Q(r))})^{ 1/3}. 
\label{9.9}
\end{align}

Appealing to \eqref{9.4}, we may choose a sequence $ \{ r_k\}$ in $ (0, R)$ such that $ r_k \rightarrow 0$ as $ k \rightarrow +\infty$, 
and 
\begin{equation}
r_k^{ -1} \| v\|^3_{ L^3(Q(r_k)} \rightarrow  0\quad  \text{ as}\quad  k \rightarrow +\infty.  
\label{9.10}
\end{equation}
By means of Jensen's inequality, having $
r_k^{ - \frac{n+2}{2}} \| v\|_{ L^2(Q(r_k))}^2 \le cr_k^{ -2/3} \| v\|_{ L^3(Q(r_k))}^2$, 
\eqref{9.10} gives  
\begin{equation}
r_k^{ -\frac{n+2}{2}} \| v\|^2_{ L^2(Q(r_k)} \rightarrow  0\quad  \text{ as}\quad  k \rightarrow +\infty.  
\label{9.11}
\end{equation}
We take $ t_k \in I(r_k)$ such that $ \| v(t_k)\|^2_{ L^2(B(r))} \le r_k^{ -\frac{n+2}{2}} \| v\|^2_{ L^2(Q(r_k))}$ for all  $ k\in \N$.  
 We now consider \eqref{9.9} with $ t=t_k, r=r_k, \eta =\eta _k$. Thanks to \eqref{9.10} and  \eqref{9.11} all terms except the first and second integral on the right-hand side of \eqref{9.9}  tend  to zero as $ k \rightarrow +\infty$. 
This shows that 
\begin{align*}
&\limsup_{k \to \infty} \intl_{ B(R)} | \widetilde{v} (t_k)|^2 \zeta dx 
\\
& \le
\limsup_{k \to \infty} \intl_{B(R)} \zeta (1- \eta_k ) d \widetilde{\sigma} _0    
= \limsup_{k \to \infty} \intl_{B(R)} \zeta (1- \eta_k ) d \widetilde{\sigma} _0.  
\end{align*}
On the other hand, by means of the weakly-$ \ast$ left  continuity  of $ \widetilde{\sigma } $, using the above inequality, we obtain  
\[
\intl_{B(R)} \zeta  d \widetilde{\sigma} _0\le \limsup_{k \to \infty} \intl_{B(R)} \zeta (1- \eta_k ) d \widetilde{\sigma} _0 = 
\intl_{B(R)} \zeta  d \widetilde{\sigma} _0 - \liminf_{ k\to \infty} \intl_{B(R)}\eta_k  d \widetilde{\sigma} _0 , 
\]
which in turn shows that 
\[
\widetilde{\sigma} _0 (\{ 0\})  \le  \lim_{k \to \infty} \intl_{B(R)} | \widetilde{v} (t_k)|^2 \zeta  dx \le  \liminf_{ k\to \infty}\intl_{B(r_k)}  \eta _k d\widetilde{\sigma} _0 = 0.
\]
Whence, the claim. \hfill \Beweisende

\subsection{Blow-up argument}

Let $\Omega \subset \Bbb R^n$. In what follows we use the following notation for the semi-norm for  the fractional derivatives of functions in the 
Sobolev-Slobodecki\u{i} spaces 
\[
| f|^p_{ W^{\theta, p}(\Omega )} = \intl_{\Omega} \intl_{\Omega}  \frac{| f(x)- f(y)|^p }{| x-y|^{ n+p\theta }} dxdy,\quad  
f\in W^{\theta, p}(\Omega ),\quad 0<\theta <1, \,\,p \ge 1. 
\]

\begin{lem}
\label{lem9.1}
Let $ v \in L^\infty(I(R); L^2 (B(R)))\cap L^3(Q(R)), 0< R<+\infty$, 
be a local suitable weak  solution to the Euler equations \eqref{9.1}.  We assume the following  local  Type I condition in terms 
of a fractional Sobolev space norm and energy concentration at time $ t=0$.
\begin{itemize}
  \item[(i)] $ \exists\,0<\theta < \frac{1}{3} $: \quad  $ \sup_{ r\in (0, R)} r^{ 3\theta - 1}  \intl_{I(r)}  | v(t)|^3_{ W^{ \theta, \, 3 }(B(r))} 
  < +\infty$.
  \item [(ii)] There exists  $ \sigma_0 \in \mathcal{M}_{ | v|^2}(0)$  with $ 
\sigma_0(\{ 0\}) >0 $. 
\end{itemize}

Then there exists a nontrivial solution $ v^{ \ast} \in L^\infty(-1, 0; L^2_\sigma (\R^n ))\cap  L^3_{ loc}([-1,0); W^{\theta ,\, 3}(\R^n ))$ to the Euler equations which fulfills the following Type I blow-up condition and energy concentration at time $ t=0$
\begin{itemize}
  \item[(iii)]  $ \sup_{ r\in (0, R)} r^{ 3\theta - 1}  \intl_{I(r)}  | v^{ \ast}(t)|^3_{ W^{ \theta, \, 3 }(B(r))} 
  < +\infty.$
  \item [(iv)] $ \mathcal{M}_{ | v^{ \ast}|^2}(0)= \{ E_0   \delta _0 \}$ for some constant $ E_0>0$. 

\end{itemize}
Furthermore, there holds the  local energy inequality for all $ \phi \in C^{\infty}_c(\R^n )$ and for 
a.e. $ -1 \le t \le s < 0$, 
\begin{align}
\intl_{ \R^n } | v^{ \ast}(t)|^2 \phi dx & \le   \intl_{ \R^n } | v^{ \ast}(s)|^2 \phi dx  	
+  \intl_{t}^{s} \intl_{\R^n }  (| v^{ \ast}|^2 + 2 p^{ \ast}) v^{ \ast}\cdot \nabla \phi  dx d\tau.  
\label{4.13}
\end{align}

\end{lem}

{\bf Proof}:  {\it 1. Scaling invariant $ L^3 $ estimate.} For notational convenience we set 
\[
K_0 = \| v(t)\|_{ L^2(-R^{ 5/2}, 0; L^2(B(R))},\quad  K_1 =  
\bigg(\sup_{ r\in (0, R)} r^{ 3\theta - 1}  \intl_{I(r)}  | v(t)|^3_{ W^{ \theta, \, 3 }(B(r))}\bigg)^{ 1/3}. 
\]
Let $ 0< r \le  R$. 
By means of H\"older's inequality together with Sobolev's inequality, we get 
\begin{align*}
\| v(s)\|^3_{ L^3(B(r))}  &\le \| v(s)\|_{ L^2(B(r))}^{ \frac{6\theta }{1+ 2\theta }} 
\| v(s)\|_{ L^{ \frac{3}{1-\theta }}(B(r))}^{ \frac{3}{1+ 2\theta }}
\\
&\le c r^{- n/2 }\| v(s)\|^3 _{ L^2(B(r))} + c\| v(s)\|_{ L^2(B(r))}^{ \frac{6\theta }{1+ 2\theta }}   
 | v(t)|_{ W^{3,\, \theta }(B(r))}^{ \frac{3}{1+ 2\theta } }.
\end{align*}
Integrating the  both sides over $ I(r)$, and using   the  H\"older's inequality, we obatin
\begin{align*}
&\intl_{I(r)}\| v(s)\|^3_{ L^3(B(r))} ds
\\
&\le c r^{ - n/2} \intl_{I(r)}\| v(s)\|^3_{ L^2(B(r))} d s  + 
c  \intl_{I(r)}\| v(s)\|_{ L^2(B(r))}^{ \frac{6\theta }{1+ 2\theta }}   
 | v(s)|_{ W^{3,\, \theta }(B(r))}^{ \frac{3}{1+ 2\theta } } ds
 \\
 &\le c r^{ - n/2} \intl_{-I(r)}\| v(s)\|^3_{ L^2(B(r))} d s 
 \\
 &\qquad + 
c r \bigg(
r^{ -5/2}\intl_{I(r)}\| v(s)\|^3_{ L^2(B(r))}  d s\bigg)^{ \frac{2\theta }{1+2\theta }} 
\bigg(r^{3\theta -1}\intl_{I(r)}| v(s)|^3_{ W^{3,\, \theta }(B(r))}d s\bigg)^{ \frac{1 }{1+2\theta }}
 \\
 & \le  cr^{ - n/2} \intl_{-I(r)}\| v(s)\|^3_{ L^2(B(r))} d s  + 
c r K_1 ^{ \frac{3 }{1+2\theta }} \bigg(
r^{ -\frac{n+2}{2}}\intl_{I(r)}\| v(s)\|^3_{ L^2(B(r))}  d s\bigg)^{ \frac{2\theta }{1+2\theta }}. 
\end{align*} 

Multiplying both sides  by $ r^{ -1}$, we  are led to 
\begin{align}
& r^{ -1} \| v\|_{ L^3(Q(r))} ^3
\cr
&\quad  \le c r^{ - \frac{n+2}{2}} \intl_{-I(r)}\| v(s)\|^3_{ L^2(B(r))} d s  + 
c  K_1 ^{ \frac{3 }{1+2\theta }} \bigg(
r^{ -\frac{n+2}{2}}\intl_{I(r)}\| v(s)\|^3_{ L^2(B(r))}  d s\bigg)^{ \frac{2\theta }{1+2\theta }}
\cr
&\quad  \le c K_0 r^{ - \frac{n+2}{2}} \| v\|_{ L^2(Q(r))}^2  + 
c K_0^{ \frac{2\theta }{1+2\theta }} K_1^{ \frac{3}{1+2\theta }} \Big(r^{ -\frac{n+2}{2}}\| v\|^2_{ L^2(Q(r))}\Big)^{ \frac{2\theta }{1+2\theta }}. 
\label{4.14}
\end{align}
Furthermore, from \eqref{4.14} we deduce that 
\begin{equation}
 \sup_{ 0<r \le R} r^{ -1} \| v\|_{ L^3(Q(r))} ^3 \le c (K_0 + K_1)^3.  
\label{4.15}
\end{equation}

{\it 2. Blow up argument}.  We assume there exists $ \sigma_0 \in \mathcal{M}_{ | v|^2}(0)$ with 
$ \sigma_0 (\{ 0\})>0$. Then from \eqref{4.14} and \eqref{4.15}  together with  Theorem\,\ref{thm9.1} we have  a positive 
constant $ \var >0$ such that 
\begin{equation}
r^{ -\frac{n+2}{2}} \| v\|_{ L^2(Q(r))}^2 \ge \var \quad \forall\,0< r \le  R. 
\label{4.16}
\end{equation}

Otherwise, \eqref{4.14} yields $ \liminf_{r \to 0} r^{ -1} \| v\|_{ L^3(Q(r))}^3 =0$, which by  Theorem\,\ref{thm9.1}   
would lead to the contradiction  $ 0 <  \sigma_0 (\{ 0\})=0$.

\hspace{0.5cm}
Now we take a decreasing sequence  $ \{ r_k\}$ in $ (0, R)$ such that $ r_k \rightarrow 0 $ as $ k \rightarrow \infty$. We define 
\[
v_k(x,t) =r_k^{ \frac{n}{2}} v(r_k x, r_k^{ \frac{n+2}{2}} t), \quad  (x,t)\in B_k\times (-1, 0), \quad  k\in \N,
\]
where
\[
B_k =  B(r_k^{ -1}R). 
\]

Clearly, $ v_k\in L^\infty(-1,0; L^2(B_k)) \cap L^3(B_k\times (-1,0))$ is a local suitable weak solution to the Euler equations in $ B_k \times (-1,0)$. Furthermore,  for every $ 0< \rho < +\infty$ the sequence $ \{ v_k\}_{ k \ge N}$ with $ r_N^{ -1} R \ge \rho $ is bounded 
in $ L^\infty(-1, 0; L^2(B(\rho )))\cap L^3(-1, 0; W^{\theta ,\, 3}(B(\rho )))$.  Thus, by means of the reflexivity and the Banach-Alaoglu theorem, 
using Cantor's diagonalization argument, eventually passing to a subsequence, we get a function    
 $ v^{ \ast} \in L^\infty(-1, 0; L^2(\R^n ))\cap L^3(-1, 0; W^{\theta ,\, 3}( \R^n ))$ with $ \nabla \cdot v^{ \ast}=0$ in $ \R^n \times (-1,0)$ in the sense of distributions such that for all $ 0< \rho < +\infty$
 \begin{align}
&  v_k \rightarrow v^{ \ast}  \quad  \text{{\it weakly-$ \ast$ in}}\quad  L^\infty(-1,0; L^2(B(\rho )))\quad \text{{\it as}}\quad  k \rightarrow +\infty,
  \label{4.17}
 \\
& v_k \rightarrow v^{ \ast}  \quad  \text{{\it weakly in }}\quad  L^3(-1, 0; W^{\theta ,\, 3}( B(\rho )))\quad  \text{{\it as}}\quad  k \rightarrow +\infty.  
 \label{4.18}
 \end{align} 
We now define,  
\[
\nabla p_{ h, k} = - E^{ \ast}_{ B_k} (v_k), \quad  \nabla p_{ 0, k} = - E^{ \ast}_{ B_k} ((v_k \cdot \nabla) v_k). 
\] 
Setting $ \widetilde{v} _k = v_k + \nabla p_{ h,k}$, we see that $ \widetilde{v} _k $ solves 
\begin{equation}
\partial _t \widetilde{v} _k + (v_k \cdot \nabla) v_k = - \nabla p_{ 0,k},\quad \text{ in}\quad  B_k\times (-1,0)
\label{4.19}
\end{equation} 
in the sense of distributions. By \eqref{2.5} having 
\begin{equation}
\| \nabla p_{ h,k}(t)\|_{ L^2(B_k)} \le c \| v_k\|_{ L^\infty(-1,0; L^2(B_k))},
\label{4.19a}
\end{equation}
 for a.e. $ t\in (-1,0)$,  and recalling that  $ \nabla p_{ h,k}$ is harmonic, we can apply the mean value property along with Jensen's inequality and \eqref{4.19a}  to find 
\[
\sup_{ B_k/2} | \nabla p_{ h,k }(t)| \le c r_k^{  \frac{n}{2}}R^{ - \frac{n}{2}} \| v_k\|_{ L^\infty(-1,0; L^2(B_k))}.
\]
Consequently, $  \nabla p_{ h,k } \rightarrow 0$ uniformly on $ B(\rho )\times (-1,0)$  as $ k \rightarrow +\infty$ 
for all $ 0<\rho <+\infty$.  
 Furthermore, employing the identity \eqref{A.1}, we see that for all $ 0< \rho < +\infty$ 
\begin{equation}
 \nabla p_{ h,k} \rightarrow 0  \quad  \text{{\it strongly in}}\quad L^\infty(-1,0; W^{1,\, 2}(B(\rho ))) \quad  \text{{\it as}}\quad  k \rightarrow +\infty.  
\label{4.20}
\end{equation}
Hence, togehter with \eqref{2.13} and \eqref{2.14} we find  for all $ 0< \rho < +\infty$
\begin{align}
&  \widetilde{v} _k \rightarrow v^{ \ast}  \quad  \text{{\it weakly-$ \ast$ in}}\quad  L^\infty(-1,0; L^2(B(\rho )))\quad   \text{{\it as}} \quad  k \rightarrow +\infty,
  \label{4.21}
 \\
& \widetilde{v} _k \rightarrow v^{ \ast}  \quad  \text{{\it weakly in }}\quad  L^3(-1, 0; W^{\theta ,\, 3}( B(\rho )))\quad  \text{{\it as}}\quad  k \rightarrow +\infty.
 \label{4.22}
 \end{align} 
On the other hand, from the estimate 
\[
\| p_{ 0, k}\|_{ L^{ 3/2}(B_k\times (-1,0))} \le \| v_k\|_{ L^3(B_k\times (-1,0))}^2,
\]
we infer that $ \{ p_{ 0, k}\}_{ k \ge N}$ is bounded in $ L^{ 3/2}(B(\rho )\times (-1,0))$.  
Thus, \eqref{4.19} shows that $ \{ \partial _t \widetilde{v} _k\}_{ k \ge N}$  is bounded in $ L^{ 3/2}(-1,0; W^{-1,\, 3/2}(B(\rho )))$. Taking into account that $ \{ v_k\}_{ k \ge N}$ is bounded in $ L^3(-1, 0; W^{\theta ,\, 3}(B(\rho )))$, we are in a position to apply the compactness lemma due to Simon\cite{sim}. This together with \eqref{4.20} yields   
\begin{equation}
 v _k\rightarrow v^{ \ast}  \quad  \text{{\it strongly  in}}\quad  L^2(B(\rho )\times (-1,0))\quad  \text{{\it as}}\quad  k \rightarrow +\infty.
\label{4.23}
\end{equation}
In particular,  for a.e. $ t\in (-1,0)$ and for all $ 0< \rho <+\infty$ it holds 
\begin{equation}
 v _k(t)\rightarrow v^{ \ast}(t)  \quad  \text{{\it strongly  in}}
 \quad  L^2(B(\rho ))\quad  \text{{\it as}}\quad  k \rightarrow +\infty.
\label{4.24}
\end{equation}

Furthermore, by means of Sobolev's embedding theorem it can be  checked easily that $ \{ v_k\}_{ k \ge  N}$ is bounded in $ L^q(B(\rho )\times (-1, 0))$ for some $ 3< q < +\infty$. Thus, \eqref{4.23} ensures that 
 \begin{equation}
 v _k\rightarrow v^{ \ast}  \quad  \text{{\it strongly  in}}\quad  L^3(B(\rho )\times (-1,0))\quad  \text{{\it as}}\quad  k \rightarrow +\infty.
\label{4.25}
\end{equation}
 Accordingly, $ v^{ \ast}$ is a weak solution to the Euler equations. Furthermore, since each element of the sequence  $\{ \widetilde{v} _k\}$ 
 satisfies the local energy inequality, 
 after letting $ k \rightarrow +\infty$, taking into account  \eqref{4.20}, \eqref{4.24} and \eqref{4.25},    we see that $ v^{ \ast}$ also fulfills the local energy inequality \eqref{4.13}.     

\hspace{0.5cm}
In addition, observing \eqref{4.16},  it holds 
\[
\| v_k\|_{ L^2(Q(1))}^2 \ge \var \quad  \forall\,k \in \N,
\] 
and thanks to \eqref{4.23} this inequality remains true for $ v^{ \ast}$, which shows that $ v^{ \ast} \neq 0$.

\hspace{0.5cm}
It now remains to check that $ v^{ \ast}$ fulfills the properties $ (iii)$ and $ (iv)$. First, by using the transformation formula of the Lebesgue integral from the definition of $ v_k$ it follows that for $ 0< \rho < 1$
\[
\rho ^{3\theta -1 } | v_k |^3_{ L^3(I(\rho ); W^{\theta ,\, 3}(B(\rho ))} =  
(r_k\rho) ^{3\theta -1 }  | v |^3_{ L^3(I(r_k\rho ); W^{\theta ,\, 3}(B(r\rho ))} \le K_1^3. 
\]
By the lower semi continuity of the semi norm $ | \cdot  |^3_{ L^3(I(\rho ); W^{\theta ,\, 3}(B(\rho ))} $ we find 
\[
\rho ^{3\theta -1 } | v^{ \ast} |^3_{ L^3(I(\rho ); W^{\theta ,\, 3}(B(\rho ))} \le K_1^3. 
\]

Now we shall  verify $ (iv)$.  Let $ k\in \N$ be fixed. From \eqref{2.10} by using the transformation formula for the Lebesgue integral, we obtain  the following  local energy inequality for $ \widetilde{v} _k$. It holds for  almost all 
 $ -r_k^{  - \frac{n+2}{2}} R^{\frac{n+2}{2}} < t < s <  0$ and for 
 all nonnegative $ \phi \in C^{\infty}(Q(r_k^{ -1} R))$ with $ \supp (\phi ) \subset B_k\times (-r_k^{ -\frac{n+2}{2}} R^{ \frac{n+2}{2}}, 0]$
 \begin{align}
&\intl_{ B_k} | \widetilde{v}_k (t)|^2 \phi   dx 
\cr
&\qquad \le   \intl_{ B_k} | \widetilde{v}_k (s)|^2 \phi dx  	
+  \intl_{t}^{s} \intl_{B_k}  |v_k |^2  v_k\cdot \nabla \phi  dx d\tau 
\cr
&\qquad  \qquad + 
 \intl_{t}^{s} \intl_{B_k}  2 p_{ 0, k}  \widetilde{v}_k \cdot \nabla \phi  dx d\tau 
 + 
 \intl_{t}^{s} \intl_{B_k}  v\cdot \nabla p_{ h,k}  v_k\cdot \nabla \phi  dx d\tau 
\cr
&\qquad \qquad +
 \intl_{t}^{s} \intl_{B_k}  v_k \otimes v_k: \nabla ^2 p_{ h, k}  \phi  dx d\tau.
\label{4.26}
\end{align}
Next, by $ \widetilde{\sigma} \in \mathcal{M}^+(B(R) )$  we denote the unique measure valued trace due to 
Lemma\,\ref{lem2.3} (cf. also Remark\,\ref{rem2.4}). Clearly, from the definition of $ \widetilde{v} _k$  the unique trace $ \widetilde{\sigma}  _k$ of $ | \widetilde{v} _k(\cdot )|^2$,  according to Lemma\,\ref{lem2.3},
is given by the relation  
\[
\intl_{ B_k}  \phi  d\widetilde{\sigma }_{k}(t)  	= 
\intl_{ B(R)}  \phi \Big(\frac{x}{r_k}\Big)  d\widetilde{\sigma }( r^{ \frac{n+2}{2}}t),\quad  \phi \in  C_c^0(B_k). 
\]
We set  $ \widetilde{\sigma} _{ 0,k}= \widetilde{\sigma }_k(0) $.  Clearly, the weakly-$ \ast$ left continuity 
of $ \widetilde{\sigma }_{ k} $ implies  
\begin{equation}
\widetilde{\sigma }_k (t)  \rightarrow \widetilde{\sigma}_{ 0,k} \quad  \text{{\it weakly-$ \ast$ in}}\quad  \mathcal{M}(B_k)
 \quad  \text{{\it as}}\quad  t \rightarrow 0^-.
\label{4.26b}
\end{equation}
 Thanks to \eqref{4.26b}   we may pass $ s \rightarrow 0$ in both sides of \eqref{4.26}. This leads to 
\begin{align}
&\intl_{ B_k} | \widetilde{v}_k (t)|^2 \phi   dx 
\cr
&\qquad \le   \intl_{B_k}  \phi  d \widetilde{\sigma } _{0, k}  	
+  \intl_{t}^{0} \intl_{B_k}  |v_k |^2  v_k\cdot \nabla \phi  dx d\tau 
\cr
&\qquad  \qquad + 
 \intl_{t}^{0} \intl_{B_k}  2 p_{ 0, k}  \widetilde{v}_k \cdot \nabla \phi  dx d\tau 
 + 
 \intl_{t}^{0} \intl_{B_k}  v\cdot \nabla p_{ h,k}  v_k\cdot \nabla \phi  dx d\tau 
\cr
&\qquad \qquad +
 \intl_{t}^{0 } \intl_{B_k}  v_k \otimes v_k: \nabla ^2 p_{ h, k}  \phi  dx d\tau,
\label{4.27}
\end{align}
Obviously, $ \| \widetilde{\sigma} _{ 0,k}\| \le \| \widetilde{\sigma}  _0 \|$ for all $ k\in \N$.  Thus, 
 by virtue of Banach-Alaoglu's theorem and Cantor's diagonalization argument we get a measure $ \sigma^*_ 0 \in \mathcal{M}^+ (\R^{n})$ 
 together with an increasing  subsequence $ \{ k_j\}$ such that for all $ 0< \rho < +\infty$
  \begin{equation}
\widetilde{ \sigma} _{0,  k_j}  \rightarrow \sigma^{ \ast}_0   \quad  \text{{\it weakly-$ \ast$ in}}\quad  \mathcal{M}(B(\rho ))
 \quad  \text{{\it as}}\quad  j \rightarrow +\infty.
 \label{4.28}
 \end{equation}
We claim that $ \sigma_0^{ \ast} = \widetilde{\sigma}_0   (\{ 0\}) \delta _0$.  Indeed, let $ \phi \in C^0_c(\R^n )$ be a nonnegative function. 
 We may choose $ 0< \rho < +\infty$ 
 such that $ \supp (\phi ) \subset  B(\rho )$.  Let $0< \var < \rho $ be arbitrarily chosen. Take $ \eta_\var  \in C^0_c(B(2\var  ))$ with 
 $ 0 \le \eta_\var  \le 1$ and $ \eta_\var (0)= 1 $. We find 
 \begin{align*}
  \intl_{ \R^n } \phi    d \sigma_0^{ \ast}  &=   \intl_{ B( 2\var  )} \phi \eta_\var     d \sigma_0^{ \ast} 
  + \intl_{ B( \rho )  \setminus \{ 0\}} \phi (1-\eta_\var )    d \sigma_0^{ \ast} 
   \end{align*} 
 From \eqref{4.28} we deduce that 
 \begin{align*}
 \intl_{ B( \rho )  \setminus \{ 0\}} \phi (1-\eta_\var )    d \sigma_0^{ \ast} 
  &= \lim_{j \to \infty} \intl_{ B(r_{ k_j}\rho )  \setminus \{ 0\}} \phi \Big(\frac{x}{r_{ k_j}}\Big)  
  \Big\{1-  \eta _\var 
  \Big(\frac{x}{r_{ k_j}}\Big) \Big\}d \sigma_0
  \\
  & \le \max \phi\, \lim_{j \to \infty} \widetilde{\sigma}_0   (B(r_{ k_j}\rho )  \setminus \{ 0\})
    \\
  & = \max \phi \,\widetilde{\sigma}  _0  \Big(\bigcap_{ j=1}^\infty  B(r_{ k_j}\rho )  \setminus \{ 0\}\Big)=0. 
    \end{align*} 
 Hence, 
 \[
\intl_{ \R^n } \phi    d \sigma_0^{ \ast} = \intl_{ B( 2\var  )} \phi \eta_\var     d \widetilde{\sigma}_0  \rightarrow \phi (0) \widetilde{\sigma}_0   (\{ 0\})
\quad  \text{ as}\quad  \var \rightarrow 0, 
\]
  which shows that $ \sigma_0^{ \ast}= a_0 \delta _0$, where $ a_0:=\widetilde{\sigma} _0  (\{ 0\})$.  
 
 \hspace{0.5cm}
  Observing \eqref{4.24}  there exists a set $ J$ of Lebesgue measure $ 0$ such that  \eqref{4.24} is satisfied 
  for all $ t\in [-1, 0]  \setminus J$ and the local energy inequalities \eqref{4.13} and \eqref{4.27} are satisfied for all $ s,t \in [-1, 0]  \setminus J$.  Taking $ t\in [-1, 0]  \setminus J$ in \eqref{4.27} with $ k_j$ in place of $ k$, and letting $ j \rightarrow \infty$,  
  we obtain the following local energy inequality   for all nonnegative $ \phi \in C^{\infty}_c(\R^n )$ and for all $ t\in [-1,0]  \setminus J$
 \begin{align}
&\intl_{ \R^n  } | v^{ \ast} (t)|^2 \phi   dx 
 \le  a_0\phi (0)  	
+  \intl_{t}^{0} \intl_{ \R^n }  |v^{ \ast} |^2  v^{ \ast}\cdot \nabla \phi  dx d\tau 
 + 
 \intl_{t}^{0} \intl_{ \R^n }  2 p^{ \ast}  v^{ \ast} \cdot \nabla \phi  dx d\tau. 
\label{4.29}
\end{align}
On the other hand, thanks to Lemma\,\ref{lem2.3} there exists a unique measure valued trace $ \sigma^{ \ast} 
\in L^\infty(-1, 0; \mathcal{M}^+(\R^{n}))$ for $ | v^{ \ast}(\cdot )|^2$,  which is weakly-$ \ast$ left continuous. Hence, in  \eqref{4.29} 
letting $ t \rightarrow 0^-$ with an appropriate choice of $ t$ in the Lebesgue set of $ | v^{ \ast}|^2$, we obtain for all nonnegative 
$ \phi \in C^{\infty}_c(\R^n )$
\[
\intl_{ \R^{n}} \phi d\sigma^{ \ast}(0) \le a_0 \delta _0.     
\]
 This shows that $ 0 \le  \sigma ^{ \ast}(0) \le  a_0 \delta _0$. Whence there exists a constant $ 0 \le  E_0 <+\infty$ such that 
 \[
\sigma^{ \ast}(0) = E_0   \delta _0.
\]
In fact, $ E_0>0$, otherwise the local energy inequality \eqref{4.13}  would imply that $ v^{ \ast} \equiv 0$. In fact, letting 
$ s \rightarrow 0^-$ in \eqref{4.13} we would obtain the inequality 
\[
\intl_{ \R^{n}} | v^{ \ast} (t)|^2 \phi  dx  \le  \intl_{0}^{t}  \intl_{ \R^{n}} v(s)\cdot \nabla \phi (| v^{ \ast}(s)|^2+ 2p^{ \ast}(s))  dx    ds 
\]
for all $ \phi \in C^{ \infty}_{ c}(\R^{n})$. Choosing an appropriate sequence of cut off function approximating $ 1$ we verify the claim.  
This completes the proof 
$ \mathcal{M}_{ | v^{ \ast}|^2}(0)= E_0   \delta _0$, the property (iv). 
 \hfill \Beweisende

 \subsection{Proof of Theorem\,\ref{thm1.4} completed}
 
 The proof will be completed by contradiction. To this end, we  assume  there exist 
$ \sigma _0 \in \mathcal{M}_{ | v|^2}(0)$ such that $ \sigma _0 (\{ x\}) >0$ for some $ x\in B(R)$. By a simple 
translation argument without loss of generality we can assume that $ x=0$. In particular, condition (ii) in Lemma\,\ref{lem9.1} is satisfied. 

\hspace{0.5cm}
In order to apply this lemma it only remains that condition (i) satisfied.    
  Let $ 0< r \le R$ be fixed.  From the Gagliardo-Nirenberg inequality we immediately get  
\[
\| v(t)\|_{ L^\infty(B(r))} \le c r^{ - \frac{n}{2}}\| v(t)\|_{ L^2(B(r))} +
 c \| v(t)\|^{ \frac{2}{n+2}}_{ L^2(B(r))} \| \nabla v(t)\|^{ \frac{n}{n+2}}_{ L^\infty(B(r))}
\]
with an absolute constant $ c>0$. 

\hspace{0.5cm}
Let $ 1 \le  s < \frac{n+2}{n} $. Taking both sides of the above inequality to the $ s$-th power,  integrate the result over $ I(r)$, and 
using the Type I blow-up condition in terms of  the velocity  gradient, we obtain 
\begin{align}
&\| v\|_{ L^s(I(r); L^\infty(B(r)))}^s
\le c r^{ - \frac{ns}{2}+ \frac{n+2}{2}} \| v\|_{ L^\infty(I(r); L^2(B(r)))}^s
\cr
&\qquad +  c\Big(\| v\|_{ L^\infty(I(r); L^2(B(r)))} + \sup (-t) \| \nabla v(t)\|_{ L^\infty(B(R))}\Big)^s
 \intl_{-r^{ \frac{n+2}{2}}}^{0} (-t)^{ -\frac{ns}{n+2}}   dt 
\cr
& \le c (K_0+ K_1)^s r^{ - \frac{sn}{2}+ \frac{n+2}{2}}
\label{4.33}
\end{align}
with $ c>0$ depending only on $ s$,  where
\[
K_0:= \| v\|_{ L^\infty(I(R); L^2(B(R)))},\quad  K_1=\sup (-t) \| \nabla v(t)\|_{ L^\infty(B(R))}. 
\]

By the standard interpolation argument  we easily get from \eqref{4.33} for every $ 1 \le s, q \le \infty$ with  
\begin{equation}
\frac{n+2}{2 s} + \frac{n}{q} > \frac{n}{2}
\label{4.34}
\end{equation}
the inequality 
\begin{equation}
\| v\|_{ L^s(I(r); L^q(B(r)))} \le c (K_0+ K_1) r^{  \frac{n+2}{2 s} + \frac{n}{q} - \frac{n}{2}}, 
\label{4.35}
\end{equation}
where $ c$ is a positive constant depending only on $ s$ and $ q$. 

\vspace{0.2cm}
Fix  $ 0< \theta < \frac{1}{3}$. We choose $ 1<p< +\infty$ such that 
\begin{equation}
p > \frac{3n\theta }{1- 3\theta }. 
\label{4.36}
\end{equation}
We set 
\[
q :=  p\frac{3-3\theta }{p-3\theta } \ge 2. 
\]
Clearly, $ 2 \le q < 3$ satisfies the relation 
\[
\frac{1}{3} =  \frac{1-\theta }{q} + \frac{\theta }{p}.
\]
Furthermore \eqref{4.34} ensures that  the following inequality holds true
\begin{equation}
\frac{(n+2)(1-3\theta )}{2(3-3\theta )} + \frac{n}{q} > \frac{n}{2}.
\label{4.37}
\end{equation} 
 Using the interpolation theorem between Sobolev-Slobodecki\u{i} spaces (cf. \cite[Theorem\,6.4.5, (7)]{bergh})  and H\"{o}lder's inequality, we get 
\begin{align*}
| v(t)|^3 _{ W^{3,\, \theta }(B(r))} &\le c r^{ -3\theta +n- \frac{3n}{q}} \| v(t)\|^3_{ L^q(B(r))} + c\| v(t)\|^{ 3-3\theta }_{ L^q(B(r))} \| \nabla v(t)\|^{3 \theta }_{ L^p(B(r))}
\\
&\le c r^{ -3\theta +n- \frac{3n}{q}}\| v(t)\|^3_{ L^q(B(r))} + cr^{ \frac{3n\theta }{p}}\| v(t)\|^{ 3-3\theta }_{ L^q(B(r))} 
\| \nabla v(t)\|^{3 \theta }_{ L^\infty(B(r))}
\\
&\le c r^{ -3\theta +n- \frac{3n}{q}}\| v(t)\|^3_{ L^q(B(r))} + c K_1^{ 3\theta }r^{ \frac{3n\theta }{p}}\| v(t)\|^{ 3-3\theta }_{ L^q(B(r))} (-t)^{ -3\theta }. 
\end{align*}
Integrating this inequality over  $ t\in I(r)$, and applying \eqref{4.35} with $ s=3$, we are lead to 
\begin{align}
| v|^3 _{ L^3(I(r); W^{3,\, \theta }(B(r)))} &\le c r^{ 1- 3\theta } (K_0+ K_1)^3
\cr
&\qquad + c K_1^{ 3\theta }r^{ \frac{3n\theta }{p}}  \intl_{I(r)} 
\| v(t)\|^{ 3-3\theta }_{ L^q(B(r))} (-t)^{ -3\theta } dt.
\label{4.38}
\end{align} 
In view of \eqref{4.37} we may choose $ \frac{3-3\theta }{1- 3\theta }< s <+\infty $ such that condition \eqref{4.34} is still fulfilled. Applying H\"older's inequality and appealing to \eqref{4.35}, we 
obtain 
\begin{align*}
& r^{ \frac{3n\theta }{p}}  \intl_{I(r)} 
\| v(t)\|^{ 3-3\theta }_{ L^q(B(r))} (-t)^{ -3\theta } dt
\\
&\le r^{ \frac{3n\theta }{p}}  \| v\|_{ L^s(I(r); L^q(B(r)))}^{ 3-3\theta } 
\bigg( \intl_{I(r)}  (-t)^{ -\frac{3s\theta }{s-3+3\theta }} dt\bigg)^{ \frac{s-3+3\theta }{s}}
\\
& \le cr^{ \frac{3n\theta }{p}}  \| v\|_{ L^s(I(r); L^q(B(r)))}^{ 3-3\theta }  r^{ \frac{n+2}{2}  \Big(1- 3\theta  - \frac{3-3\theta }{s}\Big) }
\\
& \le c (K_0+ K_1)^{ 3-3\theta } r^{ \frac{3n\theta }{p} } r^{ (3- 3\theta )\Big( \frac{n+2}{2s}+ \frac{3}{q} - \frac{n}{2}\Big)} r^{ \frac{n+2}{2}  \Big(1- 3\theta  - \frac{3-3\theta }{s}\Big) }
\\
& = c (K_0+ K_1)^{ 3-3\theta } r^{ \frac{3n\theta }{p} } r^{ (3- 3\theta )\Big( \frac{n}{q} - \frac{n}{2}\Big)} 
r^{ \frac{n+2}{2}  (1- 3\theta  ) } 
\\
&=  c (K_0+ K_1)^{ 3-3\theta } r^{ 1- 3\theta }. 
\end{align*}  
Inserting this inequality into the right-hand side of \eqref{4.38} and applying Young's inequality,  we arrive at  
\begin{equation}
| v|^3 _{ L^3(I(r); W^{3,\, \theta }(B(r)))} \le  c (K_0+ K_1)^{ 3 } r^{ 1- 3\theta },
\label{4.39}
\end{equation}
 which shows that condition (i) of Lemma\,\ref{lem9.1} is satisfied.  

\hspace{0.5cm}
Now, we are in a position to apply Lemma\,\ref{lem9.1} to  obtain a nontrivial limit 
$ v^{ \ast} \in L^\infty(-1, 0; L^2_\sigma (\R^n ))\cap  L^3_{ loc}([-1,0); W^{\theta ,\, 3}(\R^n ))$, which is a weak solution to the 
Euler equations in $ \R^{n}\times (-1,0)$ satisfying (iii) and (iv).   
On the other hand, by the assumption of the theorem $ v$ fulfills the local Type I blow up condition in terms of the velocity gradient. 
Since this Type I condition is invariant under the scaling,   the the limit function must enjoy the global Type I blow up condition in terms of the velocity gradient in $\Bbb R^n$. Since 
$ \sigma _0$ is a Dirac measure, however, by  application of Theorem\,\ref{thm3.1},   we need to have $ v^{ \ast} \equiv 0 $,  which  contradicts to the nontriviality of $ v^{ \ast}$.  \hfill \Beweisende  

\section{Proof of Corollary\,\ref{cor1.3} }
\label{sec:-4}
\setcounter{secnum}{\value{section} \setcounter{equation}{0}
\renewcommand{\theequation}{\mbox{\arabic{secnum}.\arabic{equation}}}}

 Let us consider the change of coordinates  $(x,t)\mapsto (y,\tau)$  from $ \Bbb \R^n \times (-1,0)$ to $ \Bbb \R^n \times (0, +\infty)$ by 
$$ y = \frac{x}{(-t)^{ \frac{2}{n+2} }}, \qquad  \tau = - \log (-t).
$$
Given a solution $(v,p)$ of the Euler equations,  the 
profile $ (V, P)$  in the energy conserving scale  is defined by the relation 
 \begin{align}
 v(x,t) &= \frac{1}{(-t)^{ \frac{n}{n+2}}} V
\Big(\frac{x}{(-t) ^{ \frac{2}{n+2}}}, - \log (-t)\Big),
 \label{1.2}
  \\
 p(x,t) &= \frac{1}{(-t)^{ \frac{2n}{n+2}}} P
\Big(\frac{x}{(-t) ^{ \frac{2}{n+2}}},  - \log (-t)\Big),  \quad  (x,t)\in \R^n \times (-1,0).
 \label{1.3}
  \end{align}
 We find that  the profile $(V,P)$ solves  the  following system: 
\begin{equation}
 \partial _\tau  V+ \frac{2}{n+2} y\cdot \nabla U + \frac{n}{n+2} V + (V \cdot \nabla ) V=
  - \nabla P, \qquad  \nabla \cdot  V =0\quad  
 \text{ in}\quad  \R^n . 
\label{1.4}
\end{equation}
One can also check easily that $v$ is a $ \lambda $-DSS solution of the Euler equations if and only if 
$V(\cdot, \tau)= V(\cdot, \tau+\frac{n+2}{2} \log \lambda )$ for all $\tau \in (0, +\infty)$.  Note that for a solution $  v\in L^\infty(-1,0; L^2_\sigma(\R^n ))$ to the Euler equation, satisfying \eqref{rate},  satisfies the  energy equality 
\begin{equation}
\| v(t)\|_{ L^2}^2 = \| v(-1)\|^2_{ L^2}:=E \quad  \forall\,t \in (-1, 0),
\label{1.3a}
\end{equation}
which implies also that
\begin{equation}
\| V(\tau)\|_{ L^2}^2 = \| V(0)\|^2_{ L^2}=E \quad  \forall\,\tau \in (0, +\infty).
\label{1.3b}
\end{equation}

We first show the following.

\begin{lem}
\label{lem2.1}
Let $ v\in L^\infty(-1, 0; L^2_{ \sigma }(\R^n ))\cap L^\infty_{ loc}([-1,0), W^{1,\, \infty}(\R^n ))$ be a $ \lambda $-DSS  solution to the Euler equations for some  $ 1< \lambda < +\infty$. 
Assume $ v$  satisfies  \eqref{1.3a}. Then,  for every  $0< r<+\infty$ it holds 
\begin{equation}
v \in C([-1, 0]; L^2( B(r)^c)), \quad \mbox{and}\quad \lim_{t\to 0} \|v(t)\|_{L^2 (B(r)^c)} =0.
\label{2.2}
\end{equation}

\end{lem}

{\bf Proof}: Let $0<r<+\infty$ be arbitrarily chosen. Using the transformation formula of the Lebesgue integral,  we 
calculate for $ t\in [-1,0)$
\begin{equation}
\| v(t)\|^2_{L^2(B(r )^c)} = 
\intl_{B((-t)^{ - \frac{2}{n+2}}r  )^c} | V (y, -\log (-t))|^2 dy = 
\intl_{ B((-t)^{ - \frac{2}{n+2}}r )^c} | V (y, \tau )|^2 dy, 
\end{equation}
where $ \tau = - \log (-t)$.  
Now, let $ (t_k)$ be any sequence in $ [-1, 0)$ such that $ t_k \rightarrow  0$ as $ k \rightarrow +\infty$. 
Then, $ \tau _k =  - \log (-t_k) \rightarrow  + \infty$ as $ k \rightarrow +\infty$.  On the other hand, since $ V$ is $ \lambda $-DSS 
 the profile, $ V$ satisfies  $V(\cdot , \tau)=V(\cdot, \tau+\frac{n+2}{2} \log \lambda) $ for all $\tau \in [0, +\infty)$. Accordingly, for every $ k\in \N$ there exists $ \widetilde{\tau } _k\in [0, \frac{n+2}{2}\log \lambda ]$ such that 
\[
V(\tau _k) = V(\widetilde{\tau } _k). 
\] 
Eventually, passing to a subsequence, we may assume $ \widetilde{\tau } _k \rightarrow  \tau _0 $ in $ [0, \frac{n+2}{2}\log \lambda ]$ as $ k \rightarrow  +\infty$. Thus, by using triangle inequality we obtain 
\begin{align}
\| v(t_k)\|_{ L^2(B(r )^c)} &= \| V(\tau _k)\|_{ L^2(B((-t_k)^{ - \frac{2}{n+2}}r  )^c)} 
=\| V(\widetilde{\tau } _k)\|_{ L^2(B((-t_k)^{ - \frac{2}{n+2}}r )^c)} 
\cr
& = \| V(\widetilde{\tau}  _k)\|_{ L^2(B((-t_k)^{ - \frac{2}{n+2}}r )^c)} - 
\| V(\tau _0)\|_{ L^2(B((-t_k)^{ - \frac{2}{n+2}}r )^c)} 
\cr
&\qquad + \| V(\tau _0)\|_{ L^2(B((-t_k)^{ - \frac{2}{n+2}}r)^c)} 
\cr
\label{profile}& \le  \| V(\widetilde{\tau}  _k)- V(\tau _0)\|_{ L^2(\R^n )}  +  
\| V(\tau _0)\|_{ L^2(\R^n   \setminus B((-t_k)^{ - \frac{2}{n+2}}r  ))}. 
\end{align}
To argue further we first note that   $ V$ solves the profile equation in a weak sense, namely
\begin{align*}
&\int_{\Bbb \R^n } V(y, \widetilde{\tau}  _k)\cdot \varphi (y)dy  - \int_{\Bbb \R^n } V(y,\tau_0)\cdot \varphi(y) dy 
\\
&\quad = \frac{2}{n+2}\int_{\tau_0} ^{ \widetilde{\tau}  _k}\int_{\Bbb \R^n } V (s)\cdot ( y\cdot \nabla ) \varphi (y)dyds
\\
 &\qquad- \frac{n}{n+2}\int_{\tau_0} ^{ \widetilde{\tau}  _k}\int_{\Bbb \R^n } V (s)\cdot\varphi (y)dyds
+\int_{\tau_0} ^{ \widetilde{\tau}  _k}\int_{\Bbb \R^n } V (s)\cdot ( V\cdot \nabla ) \varphi (y)dyds 
\end{align*}
 for all $\varphi \in C_c ^\infty (\Bbb \R^n )$ with $\nabla \cdot \varphi =0$,  from which, taking into account  of the fact that $ V\in L^\infty (0, +\infty; L^2_\sigma (\Bbb \R^n ))$,  we find easily that 
   $V(\widetilde{\tau}_k)\to V(\tau_0)$ weakly in $L^2(\Bbb \R^n )$ as $ k \rightarrow +\infty$. Thus the norm convergence \eqref{1.3b} together with weak convergence implies that 
the first term on the right hand side of \eqref{profile}  tends to zero as $ k \rightarrow +\infty$. Secondly, by the monotone convergence we see that also the   second term on the right hand side of \eqref{profile} tends zero as $ k \rightarrow +\infty$. Thus, 
\[
\| v(t_k)\|_{ L^2(B(r )^c)}  \rightarrow  0\quad  \text{ as}\quad  t \rightarrow 0.
\]
Since we have shown that $\|v(t_k)\|_{L^2 (B(r)^c)} \to 0$ and $v(t_k) \to 0$ weakly in $ L^2(\R^n )$ as $k\to \infty$, the conclusion \eqref{2.2} follows.
  \hfill \Beweisende  

\begin{cor}
\label{cor2.2}
Let $ 1< \lambda < +\infty$  and $ v\in L^\infty(-1, 0; L^2_{ \sigma }(\R^n ))\cap L^\infty_{ loc}([-1,0), W^{1,\, \infty}(\R^n ))$ be a $ \lambda $-DSS  solution to the Euler equations satisfying \eqref{1.3a}.  Then the energy is concentrated at $ (0,0)$, i.\,e.  \eqref{conc2} is satisfied with $ E_0=E$.  
\end{cor}

{\bf Proof}: Let $ \varphi \in C^0 (\R^n )$ and bounded. Let $ \var > 0$ be arbitrarily chosen. By the continuity of  $ \varphi $ 
we may choose  $ \delta >0$ such that $ \sup_{ x\in B(\delta )} | \varphi (x)- \varphi (0)| \le \var $.  Elementary,  
\begin{align*}
\intl_{\R^n } | v(t)|^2 \varphi  dx &= \intl_{\R^n } | v(t)|^2 (\varphi  - \varphi (0))dx + E  \varphi(0) 
\\
&= \intl_{B(\delta )^c} | v(t)|^2 (\varphi  - \varphi (0)) dx + 
\intl_{B(\delta )} | v(t)|^2 (\varphi  - \varphi (0))  dx  + E  \varphi (0). 
\end{align*}  
Thanks to \eqref{2.2} and the boundedness of $ \varphi $  we get 
\[
\limsup_{t \to 0} \intl_{\R^n } | v(t)|^2 \varphi  dx \le E  \var + E  \varphi (0),\quad \mbox{and}\quad
E  \varphi (0) \leq E  \var + \liminf_{t \to 0} \intl_{\R^n } | v(t)|^2 \varphi  dx.
\]
Therefore.
$$
\limsup_{t \to 0} \intl_{\R^n } | v(t)|^2 \varphi  dx \leq E  \varphi(0) \leq \liminf_{t \to 0} \intl_{\R^n } | v(t)|^2 \varphi  dx,
$$
which shows
$$
\lim_{t\to 0} \int_{\Bbb \R^n } |v(t)|^2 \varphi  dx =E  \varphi(0)= E  <\delta_0, \varphi>.
$$
  \hfill \Beweisende

$$\mbox{\bf Acknowledgements}$$
Chae was partially supported by NRF grants 2016R1A2B3011647, while Wolf has been supported 
supported by the NRF grand 2017R1E1A1A01074536.  
  The authors declare that they have no conflict of interest.   
  
  \appendix
\section{Some auxiliary lemmas }
\label{sec:-A}
\setcounter{secnum}{\value{section} \setcounter{equation}{0}
\renewcommand{\theequation}{\mbox{A.\arabic{equation}}}}
  
 Here we prove  fundamental properties of a harmonic function used in  the proof of the main theorem.
\begin{lem}
\label{lemA.1}
Let $ p\in L^2(\Omega )$ be a harmonic  function  on  $ \Omega $. Then for every $ \phi \in C^{\infty}_c(\Omega )$ and for all $ m\in \N$ it holds 
\begin{equation}
\intl_{\Omega}  | \nabla^m p|^2 \phi dx = \frac{1}{2^m} \intl_{\Omega} p^2 \Delta^m  \phi dx.
\label{A.1}
\end{equation} 
\end{lem}   

{\bf Proof}:  We show \eqref{A.1} by an inductive argument. First, \eqref{A.1}  with $m=1$ is clear by  the integration by parts. 
Assume \eqref{A.1} holds for  $ m\in \N$. Then we have  $ | \nabla ^{ m+1} p|^2 = \sum_{i=1}^{n} | \nabla^{ m} \partial _i p|^2  $. 
Using  the assumption  that \eqref{A.1}  holds for $ \partial _i p $ in place of $ p$, and  integrating it by parts, we infer 
\begin{align*}
\intl_{\Omega}  | \nabla^{ m+1} p|^2 \phi dx &= \sum_{i=1}^{n}\intl_{\Omega}  | \nabla^{ m} \partial _ip|^2 \phi dx
\\
&=\frac{1}{2^m} \intl_{\Omega} | \nabla p|^2 \Delta^m  \phi dx
=\frac{1}{2^{ m+1}} \intl_{\Omega} p^2 \Delta^{ m+1}  \phi dx. 
\end{align*}
 
 \begin{lem}
\label{lemA1a}  
Let $ U = \R^{n}  \setminus \overline{B(r)}, 0< r<+\infty$.  Let $ \zeta  \in C^{\infty}(U)$ denote a cut off 
function such that $ 0 \le\zeta \le 1$ in $ U$,  and $ | \nabla ^k\zeta | \le c r^{ -k}$, $ k=1, \ldots, n+1$. 
Then for every 
$ u\in L^{ \frac{2n}{n-1}}(U)$ which is harmonic in $ U$ it holds 
\begin{equation}
\| \nabla u \zeta\|_{ \infty} \le c r^{ -\frac{n+1}{2}} \|u \|_{ L^{ \frac{2n}{n-1}}(U)}.
\label{A.1a}
\end{equation}

\end{lem} 

{\bf Proof}: First, let $ x\in \R^{n}  \setminus B(2r)$. Applying the mean value property of harmonic functions and Jensen's inequality, we get 
\[
| \nabla u(x) \zeta (x)|  \le | \nabla u(x)| \le c r^{ -n-1} \intl_{B(x, r)} | u|  dx  \le c r^{ -\frac{n+1}{2}} \| u\|_{ L^{ \frac{2n}{n-1}}(U)}. 
\] 
Secondly, let $ x\in U\cap B(2r) $. By $ \eta  \in C^{\infty}_{ c}(B(4r))$ we denote a cut off function such that $ 0 \le \eta \le 1$  
in $ B(4r), \eta \equiv 1$ on $ B(2r)$ and $ | \nabla ^k \eta | \le c r^{ -k}$, $  k=1, \ldots, n+1$. Using Sobolev's embedding theorem, and applying Lemma\,\ref{lemA.1}  
with $ \phi = | \nabla ^j(\zeta \eta)|^2 $, $ j=1, \ldots, n$,  we estimate
\begin{align*}
| \nabla u(x) \zeta (x)|& \le \| \nabla u \zeta  \eta  \|_{ L^\infty(B(4r))} 
\\
&\le c \sum_{k=0}^{n} 
r^{ - \frac{n}{2} +k}\|\nabla ^k (\nabla u \zeta  \eta )\|_{ L^{2}(B(4r))} 
\\
&\le c\sum_{k=0}^{n} \sum_{j=0}^{k}
r^{ - \frac{n}{2} +k}\|\nabla ^{ k-j+1}  u  \nabla ^{ j}(\zeta \eta  )\|_{ L^{2}(B(4r))}
\\
&\le cr^{ - \frac{n+2}{2}} \| u\|_{ L^{2}(U\cap B(4r))} \le cr^{ - \frac{n+1}{2}} \| u\|_{ L^{ \frac{2n}{n-1}}(U)}.
\end{align*}
The assertion now follows from the above two estimates.  \hfill \Beweisende

\begin{lem}

\label{lemA.2}
Let $ \{ p_k\}$ be a sequence of harmonic functions in $ L^2(\Omega )$,  which converges weakly to some limit $ p$ in 
$ L^2(\Omega )$ as $ k \rightarrow +\infty$. Then $ p$ is harmonic and for every compact set $K \subset \Omega  $ and every multi index 
$ \alpha = (\alpha _1, \ldots, \alpha _n)$ it holds
\begin{equation}
 D^{ \alpha } p_k \rightarrow D^{ \alpha } p  \quad  \text{{\it uniformly on  $ K$}}\quad  \text{{\it as}}\quad  k \rightarrow +\infty. 
\label{A.3}
\end{equation}
\end{lem}

{\bf Proof}: By virtue of Weyl's lemma it is clear that $ p$ is harmonic in $ \Omega $. 
Let $ x\in \Omega $, and let  $ B(x, r) \subset \Omega $ be a ball. By the weak convergence and the mean value property 
of harmonic functions we obtain 
\[
p_k(x) = \intmw_{B(x, r)} p_k dy \rightarrow  \intmw_{B(x, r)} p dy = p(x)\quad  \text{ {\it as}}\quad k \rightarrow +\infty. 
\]   
This shows that $ p_k \rightarrow p$ pointwise as $ k \rightarrow +\infty$. In particular,  $ p_k \rightarrow p$ in $ L^2(\Omega ')$ as $ k \rightarrow +\infty$ 
for every $ \Omega ' \Subset \Omega $. Applying the identity  \eqref{A.1} for a suitable cut off function  $ \phi \ge 0$ it follows that 
$ p_k \rightarrow p$ in $ H^m(\Omega ')$ as $ k \rightarrow +\infty$ 
for every $ \Omega ' \Subset \Omega $ and for all $ m\in \N$. The uniform convergences \eqref{A.3} is now an immediate consequence 
of Sobolev's embedding theorem.  \hfill \Beweisende

\begin{lem}
\label{lemB.1}
For $ 0< R< +\infty$ define $ Q(R)= B(R) \times I(R), I(R)= (- R^{ \frac{n+2}{2}}, 0)$. 
Let $ f\in L^p(Q(R); \R^{n^2})$, $ 1< p< +\infty$. Let $ u\in L^p(Q(R))$,  solving the equation 
\begin{equation}
-\Delta u = \sum_{i,j=1}^{n} \partial _i \partial _j f_{ ij}\quad  \text{ in}\quad Q(R)
\label{B.1}
\end{equation}
in the sense of distributions. Assume for some  $ \lambda  \in (0,n) $ it holds   
\begin{equation}
\sup_{ 0< \rho < R} \rho ^{ - \lambda } \| f\|^p_{ L^p(Q(\rho ))} < +\infty. 
\label{B.2}
\end{equation}
Then there exists a constant $ c>0$ depending only on $ n,p$ and $ \lambda $ such that 
\begin{equation}
\sup_{ 0< \rho < R} \rho ^{ - \lambda   } \| u\|^p_{ L^p(Q(\rho ))} \le c \Big(R^{ -\lambda  }\| u\|^p_{ L^p(Q(R))}+ \sup_{ 0< \rho < R} \rho ^{ - \lambda } \| f\|^p_{ L^p(Q(\rho ))}\Big).  
\label{B.3}
\end{equation}
\end{lem}

{\bf Proof}:  By a routine scaling argument we may assume that $ R=1$. We extend $ f(t)$ by zero outside $ B(1)$, and denote this extension again by $ f$.   Clearly, the family of annalus $ U_j = 
B(2^{ j+1})  \setminus \overline{B(2^{ j-1})} , j\in \Z, j \le 1$, cover $ \overline{B(2)} $. By $ \{ \psi _j\}$ we denote a corresponding 
partition of unity of smooth  radial symmetric functions, such that  $ \sum_{j=-\infty}^{1} \psi _j =1$ on $ B(2)$ together with $ | \nabla \psi _j | \le c 2^{ -j}$ and $ | \nabla^2 \psi _j | \le c 2^{ -2j}$
for all $ j\in \Z, j \le 1$.  Let $ m\in \Z$, with $ m \le 0$ be arbitrarily chosen, but fixed.  We write $ u= u_1+ u_2+ u_3$, where 
\begin{align*}
u_1(x, t) &= \sum_{j=-\infty}^{m} P.V. \intl_{ \R^{n}} f(x-y, t):  \nabla ^2 N(y) \psi _j(y)dy,   
\\
u_2 (x, t) &= \sum_{j=m+1}^{1} P.V. \intl_{ \R^{n}} f(x-y, t)  : \nabla ^2 N(y) \psi _j(y)dy,
\\
u_{ 3}(x,t ) &= u(x, t)- u_1(x, t)- u_2(x,t),\quad  (x, t)\in Q(1),
\end{align*} 
where $ N$ stands for the Newton potential in $ \R^{n}$.  

\vspace{0.2cm}
Our aim will be to estimate the $L^p$ norm of $u_1, u_2$ and $ u_3$ over $ Q(2^m)$ separately.  

\vspace{0.2cm}
First, by triangle inequality we see that for   $ x\in B(2^{ m})$ and $ | x-y| \ge 2^{ m+2}$ we get $ | y| \ge 2^{ m+1}$. 
Thus, by Calder\'on-Zygmund inequality  we find for almost every $ t\in (-1, 0)$
\[
\| u_1(t)\|_{ L^p(B(2^m))}^p \le c \| f (t) \|^p_{ L^p(B(2^{ m+2}))}. 
\] 
Integration of both sides over $ I(2^m)$ with respect to time along with \eqref{B.2} yields  
\[
\| u_1\|_{ L^p(Q(2^m))}^p \le c \| f  \|^p_{ L^p(Q(2^{ m+2}))} \le 
c 2 ^{ m \lambda } \sup_{ 0< \rho < R} \rho ^{ - \lambda   } \| f\|^p_{ L^p(Q(\rho ))}. 
\] 

Next,  fix $ x\in B(2^{ m})$. It is readily seen that for all $ j \ge m+1$ it holds $ B(x, 2^{j+1 }) \subset 
B(2^{ j+1}+ 2^m) \subset B(2^{ j+2})$. Noting that  $ | k(y)| \le  c | y|^{ -n}$ it follows $ | k|\psi _j \le c 2^{ -jn}$. Accordingly,  by the aid of Jensen's inequality, and observing \eqref{B.2},   we estimate 
\begin{align*}
| u_2(x,t)| &\le c \sum_{j=m+2}^{1}\intl_{U_j}  f(x-y, t) 2^{ -jn} dx 
\le \sum_{j=m+2}^{1}
2^{ j \frac{n}{p}}\bigg(\intl_{B(2^{ j+2})} | f(y, t)|^p  dy \bigg)^{ \frac{1}{p}} 
\end{align*}
Taking the $ \esssup$ over $ x\in B(2^m)$, and taking the $ \| \cdot \|_{ L^p(I(2^m))}$ of both sides with respect to $ t$,  using Minkowski's inequality, and  observing \eqref{B.2},   we are led to 
\begin{align*}
\bigg( \intl_{I(2^m)}\| u_2(t)\|^p_{ L^\infty(B(2^m))} dt \bigg)^{ \frac{1}{p}} &
\le \bigg(\sum_{j=m+2}^{1}
2^{ j n}\intl_{Q(2^{ j+2})} | f(y, t)|^p  dy \bigg)^{ \frac{1}{p}}
\\
& \le c\Big(\sum_{j=m+2}^{1} 2^{ - j(n-\lambda) } 
\sup_{ 0< \rho < R} \rho ^{ - \lambda  }  \| f\|^p_{ L^p(Q(\rho ))}\Big)^{ \frac{1}{p}}  
\\
&\le c 2 ^{ -m  \frac{n-\lambda }{p}}\Big(\sup_{ 0< \rho < R} \rho ^{ - \lambda   } \| f\|^p_{ L^p(Q(\rho ))}\Big)^{ \frac{1}{p}}.  
\end{align*}
Consequently,  
\[
\| u_2\|^p_{ L^p(Q(2^m))} \le c 2^{ mn }  \intl_{I(2^m)}\| u_2(t)\|^p_{ L^\infty(B(2^m))} dt 
\le c 2 ^{ m \lambda } \sup_{ 0< \rho < R} \rho ^{ - \lambda  } \| f\|^p_{ L^p(Q(\rho ))}.  
\]

 In only remains to estimate $ u_3$.  By the  definition of $ u_1$ and $ u_2$, recalling that $ f(t) \equiv 0$ on $ \R^{n}  \setminus B(1)$,  we see that for almost all $ t\in (-1,0)$ and  for all $ x\in B( 1)$ it holds 
 \begin{align}
u_1(x, t) + u_2(x,t) &= \sum_{j=-\infty}^{1} P.V. \intl_{B(2)} f(x-y, t):  \nabla ^2 N(y) \psi _j(y)dy
 \cr
&= \intl_{ \R^{n}} f(x-y,t):  \nabla ^2 N(y) dy.    
 \label{B.4}
 \end{align}
 In particular, $ u_1 + u_2$ solves \eqref{B.1} in the sense of distributions. By Weyl's lemma we deduce that  $ u_3(t) = u(t)- u_1(t)-u_2(t) $ is harmonic. Thus, 
 \begin{align*}
\| u_3\|_{ L^p(Q(2^m))}^p &\le c 2^{ m n} \| u_3\|_{ L^p(I(2^m); L^\infty( B(2^m))}^p 
 \\
&\le c 2^{ m \lambda } \Big(\| u\|_{ L^p(Q(2^m))}^p + \| f\|_{ L^p(Q(1))}^p\Big).  
 \end{align*}
 Combining the estimates of $ u_1, u_2$ and $ u_3$ we get for all $ m\in \Z$, $ m \le 0$, 
\[
2^{ - m \lambda }\| u\|^p_{ L^p(Q(2^m ))} \le c \Big(\| u\|^p_{ L^p(Q(1))}+ \sup_{ 0< \rho < R} \rho ^{ - \lambda } \| f\|^p_{ L^p(Q(\rho ))}\Big).  
\]
Taking the supremum over all $ m\in \Z$, $ m \le 0$ on the left-hand side, we obtain the assertion  \eqref{B.3}.   \hfill \Beweisende 

\bibliographystyle{siam}

\end{document}